\input amstex
\documentstyle{amsppt}
\loadbold

\font\medbf=cmbx10 scaled \magstep1
\font\bigbf=cmbx10 scaled \magstep2
\magnification=1000
\TagsOnRight
\input psfig.sty
\hsize=6in
\vsize=8in
\hfuzz=3 pt
\hoffset=0.1in
\tolerance 9000
\def\L{\text{{\bf L}}}

\def\v{\vskip 1em}
\def\vs{\vskip 2em}
\def\vsk{\vskip 4em}
\def\a{\alpha}
\def\vp{\varphi}

\def\bn{\text{\bf n}}
\def\la{\big\langle}
\def\multito{{\rightsquigarrow}}

\def\ra{\big\rangle}
\def\ov{\overline}

\def\tV{\widetilde V}
\def\R{\Bbb R}
\def\vU{\varUpsilon}
\def\A{{\Cal A}}

\def\I{{\Cal I}}
\def\J{{\Cal J}}
\def\O{{\Cal O}}
\def\P{{\Cal P}}
\def\P{{\Cal P}}

\def\Om{\Omega}

\def\D{{\Cal D}}

\def\U{{\Cal U}}
\def\bU{\hbox{{\bf U}}}

\def\G{{\Cal G}}
\def\T{{\Cal T}}
\def\forallt{\hbox{for all~}}

\def\bfn{\hbox{\bf n}}
\def\ve{\varepsilon}
\def\n{\noindent}

\def\fine{\hfill $\square$}

%
\voffset=0.4in
%
%

%

%

%
%
\def\Ber{\hbox{[Be]}}
%

%
%

%
%

%

%

%
%

%
%

%

%
%

%

%

%

%

%

%
%

%
%

%
%

%
%

%
%

%
%

%
%

%

%
%

%

%

%

%

%


\def\rhead=\vbox to 2truecm{\hbox
{\ifodd\pageno\rightheadline\else\leftheadline\fi}
\hbox to 1in{\hfil}\psill}
\tenrm
\def\rightheadline{\hfil{\eightrm Nearly Time Optimal 
Patchy Feedbacks
}\hfil\folio}
\def\leftheadline{\folio\hfil{\eightrm F. Ancona and A. Bressan}\hfil}
\def\onepageout#1{\shipout\vbox
\offinterlineskip
\vbox to 2in{\rhead}
\vbox to \pageheight
\advancepageno}
\NoBlackBoxes
%
\document
\topskip=20pt
\nopagenumbers

\null
\v
\centerline{\bigbf Nearly Time Optimal Stabilizing Patchy Feedbacks}
\vskip 2cm 
\centerline{\it Fabio Ancona$^{(*)}$ and Alberto Bressan$^{(**)}$}
\vs

\parindent 40pt

\item{(*)} Dipartimento di Matematica and C.I.R.A.M., Universit\`a di Bologna,
\item{} Piazza Porta S.~Donato~5,~Bologna~40127,~Italy.
\item{} e-mail: ancona\@ciram.unibo.it
\v

\item{(*)} Dept.~of Mathematics, Penn State University, University Park,
Pa.~16802, U.S.A.
\item{} e-mail: bressan\@math.psu.edu
\vskip 1cm

\centerline{November 2005}
\vskip 1cm

\n {\bf Abstract.} We consider the time optimal stabilization
problem for a nonlinear control system $\dot x=f(x,u)$.
Let $\tau(y)$ be the minimum time needed to steer the system from the state
$y\in\R^n$ to the origin, and call $\A(T)$ the set of initial states that
can be steered to the origin in time $\tau(y)\leq T$.
Given any $\ve>0$, in this paper we construct a patchy feedback $u=U(x)$
such that every solution of $\dot x=f(x, U(x))$, $x(0)=y\in \A(T)$
reaches an $\ve$-neighborhood of the origin within time $\tau(y)+\ve$.
\vs

%
%
\vs
\vs

\noindent
{\it 1991 AMS-Subject Classification.} \
%
%
34 A,
%
%
34 D, 
%
%
49 E, 
%
%
93 D.
\vfill\eject
\vsk

\pageno=1

\parindent 20pt

\n{\medbf 1 - Introduction}
\v
Consider an optimization problem
for a nonlinear control system of the form
$$
\dot x = f(x,u)\qquad\qquad u(t)\in \bU\,,
\tag 1.1
$$
where $x\in\R^n$ describes the state of the system, the upper dot denotes
a derivative w.r.t.~time, while~$\bU\subset\R^m$ is the set of admissible
control values.
A central issue in the theory of optimal control is the
existence of a feedback control $u=U(x)$
such that  {\it all} trajectories of
$$
\dot x=f\big(x, U(x)\big)
\tag 1.2
$$
are optimal, for a given performance criterion.
In most cases,
the optimal feedback law $u=U(x)$ is not continuous.
As shown in Example~1.1 in [27] or Example~2 in [10], even
near-optimal feedback laws can usually be found only
within a class of discontinuous functions.

Therefore, it is essential to provide
suitable definitions
of ``generalized solutions'' for discontinuous ODE's.
In particular, we recall the concept of ``sample-and-hold'' solutions and Euler
solutions (limits of sample-and-hold solutions), which were
successfully implemented both within the context
of stabilization problems [16, 31, 34] and of near-optimal feedbacks
[17, 19, 27] (see also [26] for a discussion of further definitions
of generalized solutions relevant for optimization problems).
A drawback of this approach is that, as illustrated by
Example~5.3 and Example~5.4 in [30],
{\it arbitrary} discontinuous feedback
can generate too many trajectories, some of which fail to be
optimal. In fact, Example~5.3 in [30] shows that
the set of Carath\'eodory solutions
of the optimal closed-loop equation (1.2) contains,
in addition to all optimal trajectories, some other arcs that are
not optimal. Moreover, Example~5.4 in [30] exhibits an optimal control problem
in which the optimal trajectories are Euler solutions,
but the closed-loop equation (1.2) has many other Euler solutions which
are not optimal.
\v
A different strategy, proposed by Piccoli [28] and Sussmann [35],
takes as primary object of investigation an optimal ``synthesis''
which is just a collection of optimal trajectories not
necessarily arising from a feedback control.
A general notion of regular synthesis is discussed in [30]
where a sufficiency theorem for optimal synthesis is proved.
The existence and the structure of an optimal
synthesis has been the subject of a
large body of literature on nonlinear control.
At present, detailed results are
known for time optimal planar systems of the form
$$
\dot x=f(x)+g(x)\,u\qquad\qquad u\in [-1,1]\,,\qquad
x\in\R^2\,,
$$
see [9] and the references therein.
For more general classes of optimal control problems,
or in higher space dimensions,
the construction of an optimal synthesis faces severe
difficulties.
On one hand, the optimal synthesis can have an extremely
complicated structure, and only few regularity results are
presently known (see [23]).
Already for systems in two space dimensions,
an accurate description of all generic singularities of a
time optimal synthesis involves the classification
of eighteen topological equivalence classes of singular points [28, 29].
In higher dimensions, an even larger number of different
singularities arises, and the optimal synthesis can exhibit
pathological behavior such as the
the famous ``Fuller phenomenon''  (see [25], [36]), where every optimal
control has an infinite number of switchings.
On the other hand, even in cases where a regular
synthesis exists, the performance
achieved by the optimal synthesis may not be
robust.  In other words, small perturbations can
greatly affect the behavior of the synthesis (e.g. see Example~5.3 in [30]).
\v

Because of the difficulties faced in the construction of an
optimal syntheses, it seems natural to slightly relax our requirements,
and look for nearly-optimal feedbacks instead.
This is
indeed the main purpose of the present paper.
Within this wider class, one can hope to find
a feedback law whose discontinuities are sufficiently "tame",
providing the existence of trajectories in the usual
Carath\'eodory  sense, all of which
are ``almost optimal''.  Moreover, the new feedback laws
will have a simpler structure
and better robustness properties than a regular synthesis.

For sake of definiteness, we shall study the problem
of steering the system (1.1) from any initial state $y\in\R^n$ to the origin
in minimum time, under the basic assumptions
\v
\item{\bf (H)} The set $\bU\subset \R^m$ of admissible control values is bounded.
Moreover, the function $f:\R^n\times \R^m\mapsto\R^n$ is
twice continuously differentiable and has sublinear growth:
$$
\big| f(x,u)\big| \leq c\,\big(1+|x|\big)\qquad\qquad\hbox{for all}
~u\in \bU\,.
\tag 1.3
$$

\n
For $y\in \R^n$, call $T(y)$
the minimum time needed to steer the system from the state
$y\in\R^n$ to the origin, i.e. set
$$
\eqalign{ T(y)\doteq \inf \big\{ t\geq 0\,;&
~~\hbox{there exists some trajectory $x(\cdot )$ of (1.1)}\cr
&~~\hbox{that satisfies $x(0)=y,\ x(t)=0$}\big\}\,.
\cr}
\tag 1.4
$$
Roughly speaking, our main theorem states the following.
If we relax a bit the optimality requirements,
asking that every initial state $y$ be steered inside an $\ve$-neighborhood
of the origin within time $T(y)+\ve$, then this can be
accomplished by a {\bf patchy feedback}, for any fixed $\ve>0$.

\v
Patchy feedback controls
were first introduced in [1] in order to study
asymptotic
stabilization problems.  They have a particularly simple structure,
being piecewise constant in the state space $\R^n$. Moreover, the
Carath\'eodory solutions of the corresponding Cauchy problems
(1.2) enjoy important robustness properties [2, 3, 4],
which are particularly
relevant in many practical situations. Indeed, one of the main
reasons  for using a state feedback  is precisely the fact that
open loop controls are usually very sensitive to disturbances.
In particular,
we have shown in~[2]
that a patchy feedback is
``fully robust'' with respect to perturbation of the external dynamics,
and to
measurement errors having sufficiently
small total variation so to avoid the
chattering behavior that may arise at
discontinuity points.

We recall here the main definitions (see~[1]):
\v
\n{\bf Definition 1.1.} By a {\it patch} we mean a pair
$\big(\Omega,\, g\big)$ where
$\Omega\subset\R^n$
is an open domain with smooth boundary
$\partial\Om,$ and $g$ is a smooth vector field
defined on a neighborhood of the closure
$\ov\Omega$ of $\Om,$ which points strictly inward at each
boundary point $x\in\partial\Om$.
\v
Calling $\bn(x)$ the outer normal at the boundary
point $x$, we thus require
$$
\la g(x),~\bn(x)\ra <0\qquad\forallt\quad x\in\partial\Om.
\tag 1.5
$$
\v
\n{\bf Definition 1.2.}  \ We say that $g:\Omega\mapsto\R^n$
is a {\it patchy vector field} on the open domain $\Om$
if there exists a family of patches
$\big\{ (\Omega_\alpha,~g_\alpha) ;~~ \alpha\in\A\big\}$ such that

\n - $\A$ is a totally ordered set of indices,

\n - the open sets $\Omega_\alpha$ form a locally finite covering of $\Omega$,

\n - the vector field $g$ can be written in the form
$$
g(x) = g_\alpha (x)\qquad \hbox{if}\qquad x \in
\Omega_\alpha \setminus
\displaystyle{\bigcup_{\beta > \alpha} \Omega_\beta}.
\tag 1.6
$$
We shall occasionally adopt the longer notation
$\big(\Omega,\ g,\ (\Omega_\alpha,\,g_\alpha)_{_{\alpha\in \A}} \big)$
to indicate a patchy vector field, specifying both the domain
and the single patches.
\v
By setting
$$
\alpha^*(x) \doteq \max\big\{\alpha \in \A~;~~ x \in \Omega_\alpha
\big\},
\tag 1.7
$$
we can write (1.6) in the equivalent form
$$
g(x) = g_{_{\alpha^*(x)}}(x) \qquad \forallt~~x \in \Omega.
\tag 1.8
$$
\n{\bf Remark 1.1.}
Notice that the patches $(\Omega_\alpha,\,g_\alpha)$
are not uniquely determined by a patchy vector field~$(\Omega,\, g)$.
Indeed,
whenever $\alpha<\beta$, by (1.6) the values of $g_\alpha$ on the set
$\Omega_\alpha\cap\Omega_\beta$ are irrelevant.
Therefore, if the open sets
$\Omega_\alpha$
form a locally finite covering of $\Omega$ and we assume that,
for each
$\alpha\in \A$, the vector field $g_\alpha$ satisfies~(1.5)
at every point $x\in\partial\Omega_\alpha\setminus \bigcup_{\beta>\alpha}
\Omega_\beta$, then the vector field $g$
defined according with~(1.6)
is again a patchy vector field.
To see this, it suffices to construct vector fields $\tilde g_\alpha$
(defined on a neighborhood of $\overline\Omega_\alpha$ as $g_\alpha$) which
satisfy the inward pointing
property~(1.5) at every point $x\in \partial\Omega_\alpha$
and
such that $\tilde g_\alpha = g_\alpha$
on~$\Omega_\alpha\setminus \bigcup_{\beta>\alpha}
\Omega_\beta$ (cfr.[1,~Remark~2.1]).
In fact, with the same arguments one deduces that, to guarantee
that a vector field $g$ defined on an open domain~$\Omega$
according with~(1.6) be a patchy vector field,
it is sufficient to require
that each vector field $g_\alpha$ satisfy~(1.5)
at every point $x\in\partial\Omega_\alpha\setminus \big(\big(\bigcup_{\beta>\alpha}
\Omega_\beta\big)\cup \partial\Omega\big)$.

\v
If $g$ is a patchy vector field, the differential equation
$$
\dot x = g(x)
\tag 1.9
$$
has several useful properties.
In particular, in [1] it was proved that the set
of Carath\'eodory solutions of
(1.8) is closed (in the topology of uniform convergence) but possibly not
connected.
Moreover, given an initial condition
$$
x(t_0)=x_0,
\tag 1.10
$$
the Cauchy problem (1.9)-(1.10) has at least one forward
solution, and at most one backward solution, in the Carath\'eodory
sense. For every  Carath\'eodory solution \, $x=x(t)$ \, of (1.9),
the map \, $t \mapsto \a^*(x(t))$ \, is left continuous
and non-decreasing.
\v
\n{\bf Remark 1.2.}
In some situations it is useful to adopt a more general definition
of  patchy vector field than the one formulated above. Indeed,
one can consider patches $(\Omega_\alpha,~g_\alpha)$
where the domain $\Omega_\alpha$ has a piecewise smooth
boundary  (see~[3]).
In this case, the inward-pointing condition~(1.5)
can be expressed requiring that
$$
g(x)\in
\overset
\circ \to T_{\!\Omega}(x)
\tag 1.11
$$
where
$\overset \circ \to T_{\!\Omega}(x)$
denotes the interior of the
tangent  cone to  $\Omega$  at the point $x$,
defined by
$$
T_\Omega(x)\doteq
\bigg\{
v\in \R^n~:~ \liminf_{t \downarrow 0}
\frac{d\big(x+tv,\ \Omega\big)}{t}=0
\bigg\}.
\tag 1.12
$$
Clearly, at any regular point $x\in\partial\Omega$, the
interior of the
tangent cone $T_\Omega(x)$ is precisely the set of all vectors
$v\in \R^n$ that satisfy
$\big\langle v,~\bn(x)\big\rangle <0$
and hence~(1.11) coincides with the inward-pointing
condition~(1.5). One can easily see
that all the results concerning patchy
vector fields established in~[1, 2]
remain true within this more general formulation.

\v
\n{\bf Definition 1.3.}  \ Let $\big(\Omega,\ g,\
(\Omega_\alpha,\,g_\alpha)_{_{\alpha\in \A}} \big)$ be a patchy
vector field. Assume that there exist control values
$v_\alpha\in \bU$ such that, for each $\alpha\in\A,$ there holds
$$
g_\alpha(x) = f(x,\, v_\alpha)\qquad\qquad\forall~  x \in
D_\alpha \doteq
\Omega_\alpha \setminus \bigcup_{\beta > \alpha} \Omega_\beta.
\tag 1.13
$$
Then, the piecewise constant map
$$
U(x) \doteq  v_\alpha\qquad \hbox{if}\qquad x \in D_\alpha
\tag 1.14
$$
is called a \ {\it patchy feedback} \ control on $\Omega,$
and referred to as \,
$\big(\Omega,\ U,\ (\Omega_\alpha,\,v_\alpha)_{_{\alpha\in \A}} \big)$ \.
\v
\n{\bf Remark 1.3.}
By Definitions~1.2 and 1.3,
the vector field
$$
g(x)=f\big(x,\,U(x)\big)
$$
defined in connection with a given patchy feedback \
$\big(\Omega,\ U,\ (\Omega_\alpha,\,v_\alpha)_{_{\alpha\in \A}} \big)$ \
is precisely the patchy vector field \
$\big(\Omega,\ g,\ (\Omega_\alpha,\,g_\alpha)_{_{\alpha\in \A}} \big)$ \
associated with a family of fields $\big\{g_\alpha : \alpha\in \A\big\}$
satisfying (1.5)
Notice that,
recalling the notation (1.7), for all $x \in \Omega$ we have
$$
U(x) = v_{\alpha^\ast(x)}\,.
\tag 1.15
$$
As observed in Remark~1.1, the values of the vector
fields $f(x,\, v_\alpha)$ on the
set~$\Omega_\alpha\cap\Omega_\beta$ are irrelevant whenever
$\alpha<\beta$, and it is not necessary
that $f(x,\, v_\alpha)$ satisfy the inward-pointing condition~(1.5)
at the points of $\partial\Omega_\alpha\cap\big(\bigcup_{\beta>\alpha}
\Omega_\beta\big)$.
Moreover, all the properties of a patchy feedback continue to hold
even in the case where
we assume that the inward-pointing condition~(1.5)
fails to be satisfied at the points of 
$(\partial\Omega_\alpha \cap \Sigma)\setminus \bigcup_{\beta>\alpha}
\Omega_\beta$, for some
region $\Sigma$ of the boundary~$\partial\Omega$. Clearly, in this case
every Carath\'eodory trajectory of the patchy vector field $g$
can eventually reach the boundary~$\partial\Omega$
only crossing points of $\Sigma$.
\v
To state our main results, we first need to relax the
minimum time problem.
Call $\U$ the family of admissible control functions, i.e.~all
measurable functions $t\mapsto u(t)$,
$t\geq 0$, with $u(t)\in \bU$ almost everywhere.
For $y\in \R^n$ and $u\in\U$, we denote by $t\mapsto x(t;y,u)$
the solution of the Cauchy problem
$$
\dot x(t)=f\big(x(t),\, u(t)\big)\,,\qquad\qquad x(0)=y\,.
\tag 1.16
$$
The global existence and the uniqueness of this solution
are guaranteed by the assumptions (H).
Now fix $\ve>0$
arbitrarily small and define the penalization function
$$
\vp_\ve(x)\doteq
\cases
\displaystyle{|x|^2\over \ve^2-|x|^2}
\quad
&\hbox{if}
\quad |x|<\ve\,,
\\
\noalign{\smallskip}
~~\infty\quad
&\hbox{if}\quad |x|\geq\ve\,.
\endcases
\tag 1.17
$$
Consider the following $\ve$-approximate
minimization problem:
$$
\inf_{t\geq 0\,;~u\in\U} \Big\{ t+\vp_\ve\big(x(t;y,u)\big)\Big\}\,.
\tag 1.18
$$
We denote this infimum by $V(y)$, for every $y\in\R^n$, and refer
to $y\mapsto V(y)\in [0,\infty]$ as the {\it value function} for (1.18).
Observe that
$V(y)\leq T(y)$. Hence, for a fixed time $T>0$,
the set of points that can be steered to the
origin within time $T$ is contained in the
sub-level set
$$
\Lambda_T\doteq \big\{ y\in\R^n\,;~~V(y)\leq T\big\}\,.
\tag 1.19
$$
With the above notations, our main result can be stated as follows.
\v
\n{\bf Theorem 1.} {\it
Let the assumptions (H) hold and, given $\ve>0$, $T>0$, let
$\Lambda_T$ be the sub-level set defined in (1.19) in connection
with the  value function
$V$ for (1.18).
Then, there exists a patchy feedback control $u=U(x)$,
defined on a neighborhood of
%
$$
\Lambda_{T,\ve}\doteq \big\{ y\in\R^n\,;~~V(y)\leq T\,;~~|y|\geq\ve\big\},
\tag 1.20
$$
such that, for each $y\in \Lambda_{T,\ve}$, every Carath\'eodory solution
of
$$
\dot x = f\big(x,\,U(x)\big)\,,\qquad\qquad x(0)=y
\tag 1.21
$$
reaches the ball
$$B_\ve\doteq \big\{ x\in\R^n\,;~~|x|\leq\ve\big\}\,$$
within time $V(y)+\ve$.}
\v
The assumptions (H) are very general.   They do not even
imply the existence
of optimal controls, even for the relaxed problem (1.18).
We recall that the standard existence theory
requires the additional assumptions
\v
\item{\bf (H$^\prime$)} The set $\bU\subset\R^m$ of admissible
control values is compact.  For every $x\in\R^n$, the set of
velocities $\big\{ f(x,u)\,;~~u\in \bU\big\}$ is convex.
\v
If both (H) and (H$^\prime$) hold, then the infimum in (1.4) and in (1.18)
are actually attained
(e.g.~cfr.~[14]). Moreover,
the minimum time function $T:\R^n\mapsto [0,\infty]$
is lower semicontinuous.
This fact is a well known consequence of the closure property
of the graph of the set valued map $S : [0,\,\infty) \times \R^n~\multito~\R^n$
defined by $S(t,y)\doteq\big\{x(t;y,u)~;~u\in\U\big\}$\,.
Because of the lower semicontinuity of the minimum time
function, and by (1.3), it follows that,
for every $\tau\geq 0$, the
attainable set
$$
A(\tau)\doteq \big\{ y\in\R^n\,;~~T(y)\leq \tau\big\}
\tag 1.22
$$
is compact.
Since $V(y)\leq T(y)$ for all $y\in\R^n$,
from Theorem~1 one thus obtains
\v
\n {\bf Corollary.} {\it Let the assumptions (H) and (H$'$) hold, and
let $\ve>0$, $\tau>0$  be given.   Then there exists a
patchy feedback control
$u=U(x)$, defined
on a neighborhood of the
set
$$
A_\ve(\tau)\doteq \big\{ y\in\R^n\,;~~T(y)\leq \tau \,;~~|y|\geq\ve\big\}\,,
\tag 1.23
$$
%
such that, for each $y\in A_\ve(\tau)$, every Carath\'eodory solution
of (1.21)
reaches the ball $B_\ve$ within time~$T(y)+\ve$.}
\v
In
all previous papers [1, 2, 3] the construction of a stabilizing
patchy feedback
did not make any use of a control-Lyapunov
function for (1.1). Instead, the feedback law
was obtained by patching together a finite number of
open-loop controls.
We remark that a straightforward adaptation
of this strategy would not work here. Indeed, let $\ve>0$ be given.
As in [1], we can then cover the set $A_\ve(\tau)$ with finitely many
tubes $\Omega_1,\ldots,\Omega_N$ and construct a patchy feedback
$u=U_\alpha(x)$
steering each point $y\in\Omega_\alpha$ inside the ball $B_\ve$
within time $T(y)+\ve$.
However, we cannot guarantee that the patchy feedback
$$u(x)=U_{\alpha^*(x)}\qquad\qquad \alpha^*(x)\doteq \max\,\{\alpha\,;~~
x\in \Omega_\alpha\}
\tag 1.24
$$
is nearly-optimal (see Fig.1).
Indeed, call $T_\alpha(y)$ the time taken by the control $U_\alpha$
to steer the point $y\in\Omega_\alpha$ inside $B_\ve$.
Let $t\mapsto x(t)$ be a trajectory of the patchy feedback (1.24),
with $x(0)=y$, $x(\tau)\in B_\ve$.
Assume $\alpha^*(t)=\alpha$ for $t\in \,]t_{\alpha-1}\,,\,t_\alpha]\,$.
{}The near-optimality of each feedback implies
$T_\alpha(x)\leq T(x)+\ve\,$  for every $x$.  Moreover
$$T_\alpha(x(t_{\alpha-1}))-T_\alpha(x(t_\alpha)) = (t_\alpha-t_{\alpha-1})\,.$$
Unfortunately, from the above inequalities one can only deduce
$$T(x(t_{\alpha-1}))-T(x(t_\alpha)) \geq (t_\alpha-t_{\alpha-1})-\ve\,.$$
and hence
$\tau\leq T(y)+N\ve\,$. This is a useless information, because
the number $N$ of tubes may well approach infinity as $\ve\to 0$.

\midinsert
\vskip 10pt
\centerline{\hbox{\psfig{figure=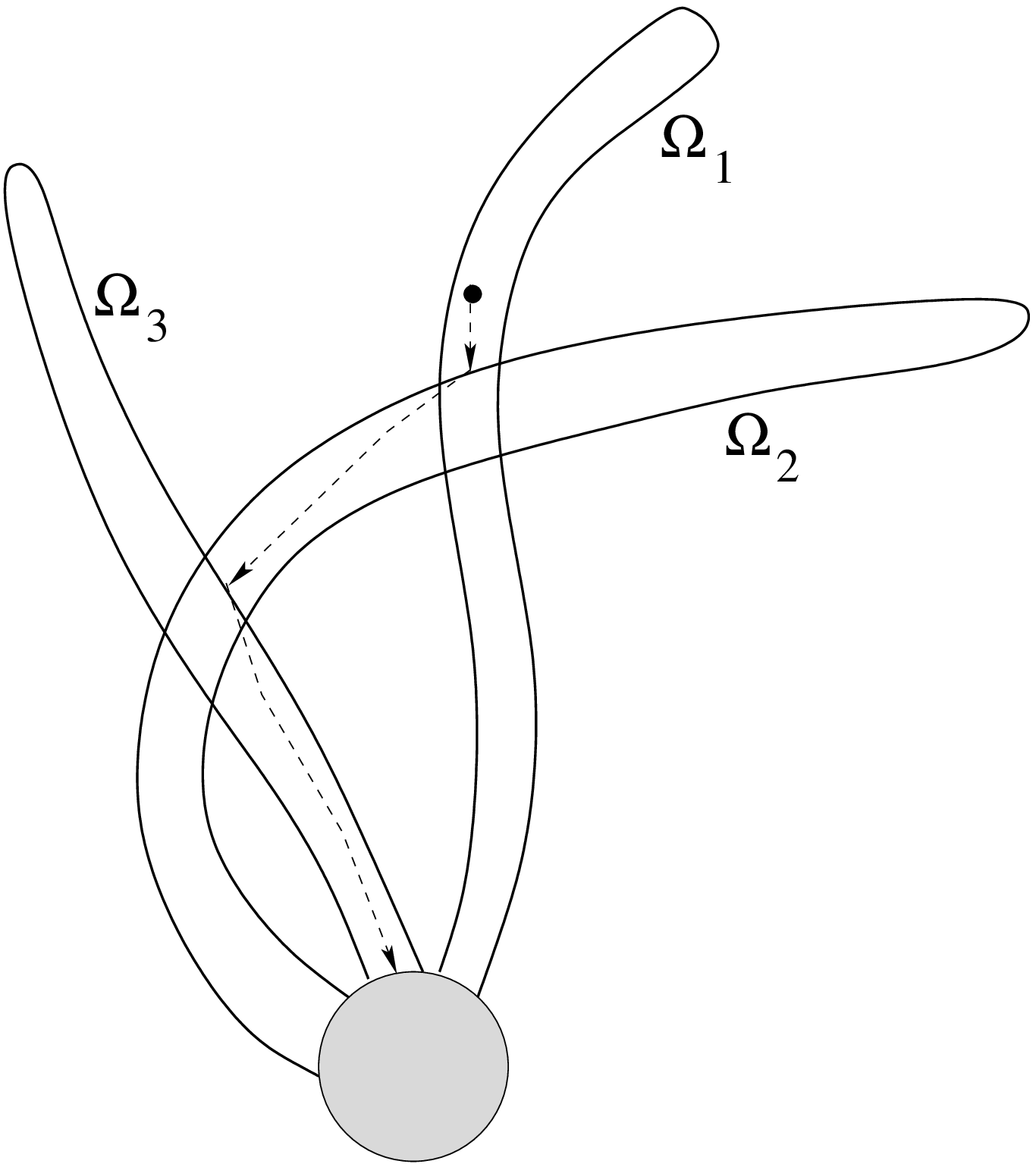,width=7cm}}}
\centerline{\hbox{\smc figure 1}}
\vskip 5pt
\endinsert
\n

To overcome this problem, in the present paper we perform
an entirely different
construction of the patchy feedback.  As starting point, instead of
open-loop controls, we use the
value function $V$ for the problem
(1.18), together with a piecewise quadratic approximation $\widetilde V$.
This has the form
$$\widetilde V(x)=\min_j V_j(x)\,,\qquad\qquad
V_j(x)=a_j+b_j\cdot x+ c\,|x|^2.$$
and satisfies $\widetilde V(x)\leq V(x)+\ve$ for each point $x$.
The result will be achieved by constructing a patchy feedback
such that
$${d\over dt} \widetilde V \big(x(t)\big)=
\nabla \widetilde V\big(x(t)\big)\cdot f\big(x(t),\,u(x(t)\big)
\leq \ve$$
at a.e.~time $t$.

\vsk
\n{\medbf 2 - Preliminary results}
\v
Throughout the paper,
by $B(x,r)$ we denote the closed ball centered at
$x$ with radius $r$, and set $B_r\doteq B(0,r)$. The closure, the interior and the
boundary of a set
$\Omega$ are written as $\overline\Omega$,
$\overset \circ \to \Omega$ and~$\partial\Omega$, respectively,
while $\text{\it diam}(\Omega)$ denotes the diameter of a bounded set $\Omega$.
The distance of a point $x$ from a set $\Omega$ is denoted by $d_\Omega(x)$,
while $d_\Omega(E)\doteq\inf_{x\in E}d_\Omega(x)$ denotes the distance between two sets~$\Omega,\, E$.
The number of elements of a finite set $\J$ is denoted by $|\J|$.

We begin by observing that the infimum in (1.18)
provides an upper bound for the time
needed to steer the system (1.1) from $y$ to the ball $B_\ve$.
Hence, for every $T\geq 0$, the sub-level set $\Lambda_T$ 
of the value function $V$ for (1.18)
is contained in the set of points that can be steered to the 
ball $B_\ve$ within time $T$.
On the other hand, notice that 
the scalar Cauchy problem
$$
\dot z = c\,(1+z)\,,\qquad\qquad z(0)=\ve\,,
\tag 2.1
$$
has solution
$$
z(t)=   (1+\ve)e^{c t}-1\,.
\tag 2.2
$$
Therefore, because of  (1.3),  a comparison argument yields
$$
\Lambda_T\subseteq B_{z(T)}\,,
\tag 2.3
$$
for every $T\geq 0$.
\v
In connection with the relaxed minimization problem (1.18), we now 
show that  the value function $V$ is
Lipschitz continuous  on $\Lambda_T$ and locally semiconcave,
that is, for any $x_0$,
there exists a constant~$c_0>0$ such that
there holds
$$
V(x_1)+V(x_2)-2 V\bigg(\frac{x_1+x_2}{2}\bigg)\leq c_0\big|x_1-x_2\big|^2
\tag 2.4
$$
for all $x_1,\, x_2$ in a neighborhood of $x_0$.
We refer to [14] for the definition and properties of semiconcave
functions.
\v
\n{\bf Lemma 1.} {\it With the assumptions (H),
for any fixed $\ve,\,T>0$ the restriction of
the value function~$V$ for~(1.18) to the sublevel set $\Lambda_T$ is
Lipschitz continuous and locally semiconcave.
Indeed, there exists a positive constant $\lambda$ such that,
for every point $y_0\in \Lambda_T$ where $V$ is differentiable, there holds}
%
$$
V(y)\leq V(y_0)+\la \nabla V(y_0),\, y-y_0\ra+\lambda\, |y-y_0|^2\qquad
\qquad \forall~y\in\Lambda_T\,.
\tag 2.5
$$
\v
\n{\bf Proof.} 

\n{\bf 1.} First observe that, since we are only proving something about the
value function $V$ for (1.18), it is not restrictive
to assume that the additional hypotheses (H$'$) hold.
Indeed, allowing the set of controls to range in
the closure of $\bU$ does not affect the value function.
Moreover, if the sets of velocites 
$\big\{ f(x,u)\,;~~u\in \bU\big\}$ are not convex,
we can replace the original system (1.1) by a chattering one (see~\Ber ),
such that the problem (1.18) yields exactly the same value function.
This in particular implies that the value function $V$ is lower 
semicontinuous and that the sub-level set (1.19) is compact.

\v
\n{\bf 2.} Next, observe that, since the function $f$ is twice continuously differentiable
and the sets $\Lambda_T$, $\bU$ are compact,
by standard differentiability properties of the
trajectories of a control system (1.1), 
there holds
$$
\align
\sup
\Sb
t\in[0,T],\ y\in \Lambda_T
\\
\noalign{\smallskip}
u\in\U
\endSb
\bigg| {\partial\over\partial y}\, x(t;y,u) \bigg|
&\leq  \exp\big\{ T\,\|D_x f\|_{\L^\infty(\Lambda_T)}
\big\}\doteq M_1\,,
\tag 2.6
\\
\noalign{\medskip}
\sup
\Sb
t\in[0,T],\ y\in \Lambda_T
\\
\noalign{\smallskip}
u\in\U
\endSb
\bigg| {\partial^2\over\partial y^2}\, x(t;y,u) \bigg|
&\leq \Big[M_1\big(1+T M_1^3 \,\|D^2_x f\|_{\L^\infty(\Lambda_T)}\big)\Big]
\doteq M_2\,,
\tag 2.7
\endalign
$$
where
$$
\aligned
\|D_x f\|_{\L^\infty(\Lambda_T)}&\doteq\sup_{x\in\Lambda_T, ~u\in\bU} 
\bigg| {\partial\over\partial x}f(x,u)\bigg|<\infty\,,
\\
\noalign{\medskip}
\|D^2_x f\|_{\L^\infty(\Lambda_T)}&\doteq\sup_{x\in\Lambda_T, ~u\in\bU} 
\bigg| {\partial^2\over\partial x^2}f(x,u)\bigg|<\infty\,,
\endaligned
\tag 2.8
$$
provide a bound on the first and second partial derivatives
of $f$ w.r.t.~the $x$-variable,
over the set~$\Lambda_T$.
Then, because of (2.6),
there exists a constant $c_1$ such that
$$
\big| x(t;y_2,u)- x(t;y_1,u)\big|
\leq c_1|y_2-y_1|\qquad\quad\forall~t\in[0,T],\ \ y_1, y_2\in \Lambda_T,\ \ u\in\U\,.
\tag 2.9
$$
\v

\n{\bf 3.} Given $y_0\in \Lambda_T$, by the previous assumptions at point {\bf 1} there
exists an optimal control $u_0\in\U$,
and a time $t_0$,
such that
$$ 
t_0 + \vp_\ve\big(x(t_0; y_0, u_0)\big)= V(y_0)\leq T\,.
\tag 2.10
$$
This, by definition (1.17) of $\vp_\ve$\,, of course implies
$$
t_0\leq T\,,\qquad \qquad \big|x(t_0; y_0, u_0)
\big|\leq \ve\,\sqrt{{T\over 1+T}}\,.
\tag 2.11
$$
Hence, using (2.11) together with (2.9), we find that there exists some 
constant $\delta>0$, depending only on $\ve, T,$ and on $c_1$, 
but not on the point $y_0\in\Lambda_T$, such that
$$
\big| x(t_0;y,u_0)\big|\leq \ve\,\sqrt{{{2T+1}\over {2(1+T)}}}
\qquad\quad\forall~y\in B(y_0,\delta)\cap  \Lambda_T\,.
\tag 2.12
$$
Observe now that,
since $V(y)$ is the infimum in (1.18), there holds
$$
V(y)\leq V^0(y)\doteq t_0+ \vp_\ve\big(x(t_0; y, u_0)\big)
\qquad\forall~y\,.
\tag 2.13
$$
Because of (2.12), the map $y\mapsto V^0(y)$ defined in (2.13) 
is twice continuously differentiable
at every point of $B(y_0,\delta)\cap  \Lambda_T$.
Hence, since (2.10) implies $V^0(y_0)=V(y_0)$, there holds
$$
V^0(y)\leq V(y_0)+ \la\nabla V^0(y_0),\, y-y_0\ra+\lambda_0
\big|y-y_0\big|^2
\qquad\quad\forall~y\in B(y_0,\delta)\cap  \Lambda_T\,,
\tag 2.14
$$
%
with
$$
\lambda_0\doteq\sup_{\eta \in B(y_0,\delta)\cap  \Lambda_T}
\big|D^2  V^0(\eta)\big|\,.
$$
The gradient of the function $V^0$ is computed by
$$\nabla V^0(y)=
\nabla \vp_\ve\big(x(t_0;y,u_0)\big)\cdot {\partial\over\partial y}\, 
x(t_0;y,u_0)\,.
$$
Thus, relying on (2.6), (2.12), and setting 
$$
M_0\doteq \sup\Sb
t\in[0,T],\ y\in \Lambda_T
\\
\noalign{\smallskip}
u\in\U
\endSb
 \big|x(t;y,u)\big|\,,
\tag 2.15
$$
we obtain
$$
\qquad
\big|\nabla V^0(y)\big|\leq {2\ve^2\big|x(t_0;y,u_0)\big|
\over \Big(\ve^2-\big|x(t_0;y,u_0)\big|^2\Big)^2}
\cdot M_1
\leq 
{{8(1+T)^2 M_0 M_1}\over{\ve^2}}
\doteq c_2
\qquad\ \forall~y\in B(y_0,\delta)\cap  \Lambda_T\,.
\tag 2.16
$$
With similar computations, using (2.7), (2.12),
we find that a bound on the
second derivative of $V^0$ is provided by
$$
\aligned
\big|D^2 V^0(y)\big|&\leq 
\big|D^2\vp_\ve\big(x(t_0;y,u_0)\big)\big|
\cdot \bigg| {\partial\over\partial y}\, x(t;y,u) \bigg|
+ \big|\nabla \vp_\ve\big(x(t_0;y,u_0)\big)\big|\cdot 
\bigg| {\partial^2\over\partial y^2}\, x(t;y,u) \bigg|
\\
\noalign{\smallskip}
\leq
\bigg[&{{8(1+T)^2 M_1}\over{\ve^2}}+
{{64(1+T)^3 M_0^2 M_1}\over{\ve^4}}
+{{8(1+T)^2 M_0 M_2}\over{\ve^2}}
\bigg]\doteq c_3
\qquad\quad\forall~y\in B(y_0,\delta)\cap  \Lambda_T\,.
\endaligned
\tag 2.17
$$
Notice that the constants $c_2, c_3$ depend
only on $\ve, T,$ and on the function $f$, 
but not on the point $y_0\in\Lambda_T$.
Then, (2.13), (2.14), together with (2.16), (2.17) yield
$$
V(y)\leq V(y_0)+\big(c_2+\delta\,c_3)
\big|y-y_0\big|
\qquad\forall~y\in B(y_0,\delta)\cap  \Lambda_T\,,
\quad\forall~y_0\in\Lambda_T\,,
\tag 2.18
$$
which, in turn, implies 
$$
~y_1,y_2\in \Lambda_T\,,\,
\big|y_1-y_2|<\delta\quad\Longrightarrow\quad
\big|V(y_1)-V(y_2)\big|\leq \big(c_2+\delta\,c_3)\cdot\big|y_1-y_2|\,.
\tag 2.19
$$
Since the set $\Lambda_T$ is compact, we deduce from (2.19) 
that the map $V$ is (globally)
Lipschitz continuous on $\Lambda_T$.
\v

\n{\bf 4.} Given $y_0\in \Lambda_T$, in connection
with the constants $\lambda_0, \delta, c_2$ introduced  at point {\bf 3}
choose 
$$
\lambda>\max\bigg\{\lambda_0,\, \frac{\text{Lip}(V)+c_2}{\delta}\bigg\}\,,
$$
and observe that, because of (2.16), there holds
$$
\aligned
\la\nabla V^0(y_0),\, y-y_0\ra+\lambda\big|y-y_0\big|^2
&\geq (-c_2+\lambda\,\delta)\big|y-y_0\big|
\\
\noalign{\smallskip}
&\geq \text{Lip}(V)\cdot\big|y-y_0\big|
\endaligned
\qquad\quad
\forall~y\in\Lambda_T\setminus B(y_0,\delta)\,.
\tag 2.20
$$
Thus, (2.13), (2.14), together with (2.20), yield
$$
V(y)\leq V(y_0)+ \la\nabla V^0(y_0),\, y-y_0\ra+\lambda\,
\big|y-y_0\big|^2
\qquad\quad\forall~y\in \Lambda_T\,.
\tag 2.21
$$
\v

\n{\bf 5.} Fix $\rho>0$. By the above arguments there exist a positive
constant
$\lambda=\lambda_\rho$ so that, for every fixed $y_0\in \Lambda_{T+\rho}$,
the estimate (2.21)
holds for all $y\in \Lambda_{T+\rho}$.
Next, given $x_0\in \Lambda_T$, choose $\delta_0$ so that
$B(x_0,\delta_0)\subset \Lambda_{T+\rho}$.
Then, for any $y_1,\, y_2 \in B(x_0,\delta_0)\cap  \Lambda_T$,
one has $\frac{y_1+y_2}{2}\in \Lambda_{T+\rho}$.
Hence, applying~(2.21)  for $y=y_i$, $i=1,2$, and $y_0=\frac{y_1+y_2}{2}$,
and summing up the corresponding inequalities,
since $y_1-y_0=y_0-y_2$ we obtain
$$
V(y_1)+V(y_2)\leq 2 V\bigg(\frac{y_1+y_2}{2}\bigg)+ 
\lambda\Big[\big|y_1-y_0\big|^2+\big|y_2-y_0\big|^2\Big]
\leq 2 V\bigg(\frac{y_1+y_2}{2}\bigg)+ \lambda\big|y_1-y_2\big|^2\,,
$$
which shows that the estimate~(2.4) is verified, with $c_0=\lambda$,
for all $y_1,y_2\in \Lambda_T$
in the ball $B(x_0,\delta_0)$.
Therefore, the map $V$ is locally semiconcave on $\Lambda_T$.
\v

\n{\bf 6.} 
To conclude the proof of the lemma, consider a point $y_0\in \Lambda_T$
where $V$ is differentiable, and observe that, by (2.21), one has
$$
\Big\langle
\nabla V(y_0)-\nabla V^0(y_0),\
\displaystyle{\frac{y-y_0}{|y-y_0|}}\Big\rangle\leq
\lambda\bigg[|y-y_0|+\displaystyle{\frac{o(|y-y_0|)}{|y-y_0|}}
\bigg]
\tag 2.22
$$
for all $y\in \Lambda_T$.
Thus, taking $y_\sigma\doteq y_0 +\sigma\, \frac{\nabla V(y_0)-\nabla V^0(y_0)}
{|\nabla V(y_0)-\nabla V^0(y_0)|}$, \, $\sigma>0$
from (2.22) we deduce
$$
\big|\nabla V(y_0)-\nabla V^0(y_0)\big|\leq
\lambda\Big[\sigma+\displaystyle{\frac{o(\sigma)}{\sigma}}\Big]
\qquad \forall~\sigma>0\,.
\tag 2.23
$$
By letting $\sigma\to 0$ in (2.23) we obtain $\nabla V^0(y_0)=\nabla V(y_0)$
which, together with (2.21), yields~(2.5), completing the proof
of the lemma.
\fine
\vs

We next show that the  value function $V$ enjoys an infinitesimal decrease property
at every point where it is differentiable, which is expressed
in terms of an
Hamilton-Jacobi inequality.
\v

\n{\bf Lemma 2.} {\it With the assumptions (H),
given $\ve,\,T>0$, let $V$ be the  value function for (1.18). 
Then, there exists $0<\ve_0<\ve$ such that, letting $\Lambda_{T,\ve_0}$
be the set defined in (1.20), for each $y\in\Lambda_{T,\ve_0}$ at which $V$ is differentiable
there holds
$$
\inf_{v\in\bU}\Big\{\big\langle\nabla V(y),\, f(y, v)\big\rangle\Big\}+1\leq 0\,.
\tag 2.24
$$
}
\v
\n{\bf Proof.} 
Given $\ve,\,T>0$, set
$$
\ve_0\doteq \ve\,\sqrt{\frac{4T+1}{2(1+2T)}}\,,\qquad\qquad
\ve'_0\doteq\ve\,\sqrt{{2T\over 1+2T}}\,,
\tag 2.25
$$
$$
\tau_0\doteq c^{-1}\ln\bigg(\frac{\ve_0+1}{\ve'_0+1}\bigg)\,,
\tag 2.26
$$
where $c$ denotes the constant in (1.3),
and observe that, by definition (1.17) of $\vp_\ve$\,,
one has
$$
\vp_\ve(x)\geq 2T
\qquad\quad\text{whenever}\quad\qquad |x|\geq\ve'_0\,.
\tag 2.27
$$
Then, recalling that (2.2) provides the solution to the scalar
Cauchy problem (2.1), by a comparison argument, and because of (1.3),
we deduce that
$$
|y|\geq\ve_0\qquad\Longrightarrow\qquad \big|x(t;y,u)\big|\geq\ve'_0
\qquad\forall~t\in[0,\tau_0],\ \ u\in\U\,.
\tag 2.28
$$
Hence, (2.27) together with (2.28), yields
$$
t+\vp_\ve\big(x(t; y, u)\big)\geq 2T\qquad\quad\forall~t\in[0,\tau_0], \ \ |y|\geq\ve_0,
\ \ u\in\U\,.
\tag 2.29
$$
>From (2.29) we deduce that, for every $y\in\Lambda_{T,\ve_0}$,
the value function for (1.18) satisfies
$$
V(y)=\inf_{t> \tau_0\,;~u\in\U} \Big\{
t+\vp_\ve\big(x(t;y,u)\big)\Big\}> \tau_0\,.
\tag 2.30
$$
Thus, we reach the conclusion of the Lemma observing that by
standard arguments in control theory (e.g. see [14])
one can show that
the value function for (2.30) satisfies the Hamilton-Jacoby
inequality~(2.24)
at every point where $V$ is differentiable.
\fine
\v

\n{\bf Remark~2.1.} Notice that, in the proof of Theorem~1, we shall
only need to have at a disposal a value function $V$ satisfying
the conclusions of Lemma~1 and Lemma~2. 
\vs

We state now two technical results which will be useful later
in the 
construction of an almost time optimal patchy feedback.
We shall provide a proof of them in the Appendix at the end of the paper.
Throughout the following,
for any given subset $C$ of a sphere $S$, we let $\partial_{\!_{S}} C$
denote the boundary of $C$ relative to the topology
of $S$.
\v

\n{\bf Lemma 3.} 
{\it Given $r_0>0$, let $S$ be a sphere with radius $r\geq r_0$, and let
$g$ be a bounded,\linebreak
Lipschitz continuous vector field which points strictly
inward at  the points of a closed 
set  $C\subset S$ that has a piecewise smooth
relative boundary $\partial_{\!_{S}} C$.  More precisely, letting 
$\bfn_{S}(y)$ denote the unit outer normal to  $S$ at the point 
$y$, 
assume that
$$
\big\langle
\bfn_{S}(y),\,g(y)\big\rangle\leq - \overline c
\qquad\quad\forall~y\in C\,,
\tag 2.31
$$
for some constant $\overline c>0$.
Denote by $t\mapsto x(t,y)$
the solution of the Cauchy problem $\dot x= g(x)$, $x(0)=y$.
Then there exists  $\overline\ve>0$, depending only on $r_0, \overline c$,
$\|g\|_{\L^\infty}$, 
and on the Lipschitz constant $\text{Lip}(g)$ of $g$,  such that the following holds.
Define 
$$
\Gamma_{\overline\ve}(C)\doteq
\big\{\, x(\tau,y)\,; ~~y\in B(C,\,\overline \ve)\cap S\,,~~
{d_{\!_{C}}}^{\!\!2}(y)-\overline\ve^2
< \tau\leq 0\,\big\}\,.
\tag 2.32
$$
Then the vector field $g$ is transversal to the boundary of $\Gamma_{\overline\ve}\doteq
\Gamma_{\overline\ve}(C)$.
Indeed, it points strictly inward on the set
$$
\partial^-\Gamma_{\overline\ve}\doteq
\big\{\, x({d_{\!_{C}}}^{\!\!2}(y)-\overline\ve^2,\,y)\,; ~~y\in \overset\ \circ\to
B(C,\,\overline \ve)\cap
S
\,\big\}
\tag 2.33
$$
and strictly outward on the set
$$
\partial^+\,\Gamma_{\overline\ve}\doteq
\partial\,\Gamma_{\overline\ve}\cap S\,.
\tag 2.34
$$
}

\v

The lens-shaped domain (2.32) provides the basic building block
for the construction of the patchy feedback produced in the next section. 
In some situations 
it will be necessary to restrict such domains cutting them along 
hyperplanes in
order to preserve the (almost) time-optimality property of the feedback law.
The next lemma provides an a-priori lower bound on the
distance between the upper boundary of a collection 
of such domains and the union of spheres around which 
the domains are cosntructed.
%
\midinsert
\centerline{\hbox{\psfig{figure=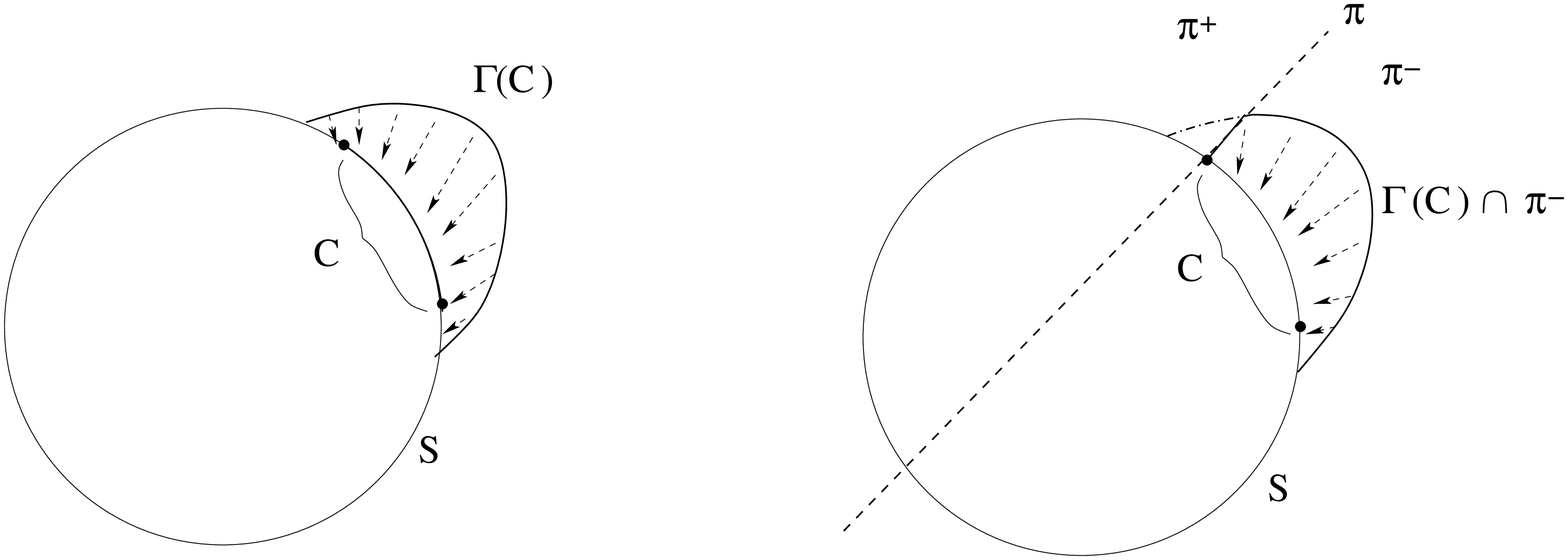,width=13cm}
}}
\centerline{\hbox{\smc figure 2}}
\vskip 5pt
\endinsert
\n

\vs

\n{\bf Lemma~4.}
{\it Given $0<r_0<r'_0$, let $B_1,\dots, B_\nu$
be a finite collection of balls with sur\-fa\-ces~$S_1,\dots, S_\nu$, having
radii $r_1,\dots, r_\nu\in [r_0,\,r_0']$,
and satisfying
$$
\aligned
S_i&\cap
\bigcup_{j=1}^\nu 
\overset\circ\to B_j
\neq \emptyset\,,
\\
\noalign{\smallskip}
S_i&\setminus\bigcup_{j=1}^\nu 
\overset\circ\to B_j
\neq \emptyset\,,
\endaligned
\qquad\quad\forall~i=1,\dots,\nu\,.
\tag 2.35
$$
Consider the sets 
$$
C_i\doteq S_i\setminus\Big(\bigcup_{j=1}^\nu
\overset\circ\to B_j\Big)
\qquad\quad \forall~i=1,\dots,\nu\,,
\tag 2.36
$$
and let $g_1,\dots, g_\nu$ be bounded,
Lipschitz continuous vector fields which point strictly
inward (towards the interior of $S_1,\dots, S_\nu$) on
$C_1,\dots, C_\nu$, respectively.
Then, there exist constants  $\overline\ve',\, c_4>0$, 
depending only on $\nu$, \, $r_0, r'_0$,
and on  $\|g_i\|_{\L^\infty}$, 
$\text{Lip}(g_i)$, $i=1,\dots,\nu$, such that the following holds.
Let 
%
$$
\gathered
\Pi=\bigcup_{k=1}^\nu\Pi_k\,,
\qquad\qquad
\Pi_k\doteq
\big\{
\pi_{k,i}\,;~~~ i\in\J_k
\big\}\,,
\\
\noalign{\medskip}
\J_k\subset \{1,\dots,\nu\}\setminus\{k\}
\qquad\forall~k\,,
\endgathered
\tag 2.37
$$
be  a (possibly empty) collection of hyperplanes
enjoying the properties:
\v
\item{--\ } 
$\pi_{k,i}=\pi_{i,k}$ \, for all \, $i\in\J_k$, $k\in\J_i$\,;

\item{--\ } 
$\pi_{k,i}\in \Pi_h$ \, if and only if either \, $h=k$, $i\in\J_h$, \, or
\, $h=i$, $k\in\J_h$\,;

\item{--\ } if $\overset\circ\to B_k\cap \overset\circ\to B_i\neq\emptyset$\,,
$i\in\J_k$,
then $\pi_{k,i}$ is the hyperplane passing through $S_k\cap S_i$\,
(cfr. Fig.~3);

\item{--\ } if $\overset\circ\to B_k\cap \overset\circ\to B_i=\emptyset$\,,
$i\in\J_k$,
then $\pi_{k,i}$ is an hyperplane separating $S_k$ and $S_i$, i.e. s.t. 
$S_k$, $S_i$ are entirely contained in the opposite closed half spaces determined
by $\pi_{k,i}$ (cfr. Fig~4)\,.
\v
\midinsert
\centerline{\hbox{\psfig{figure=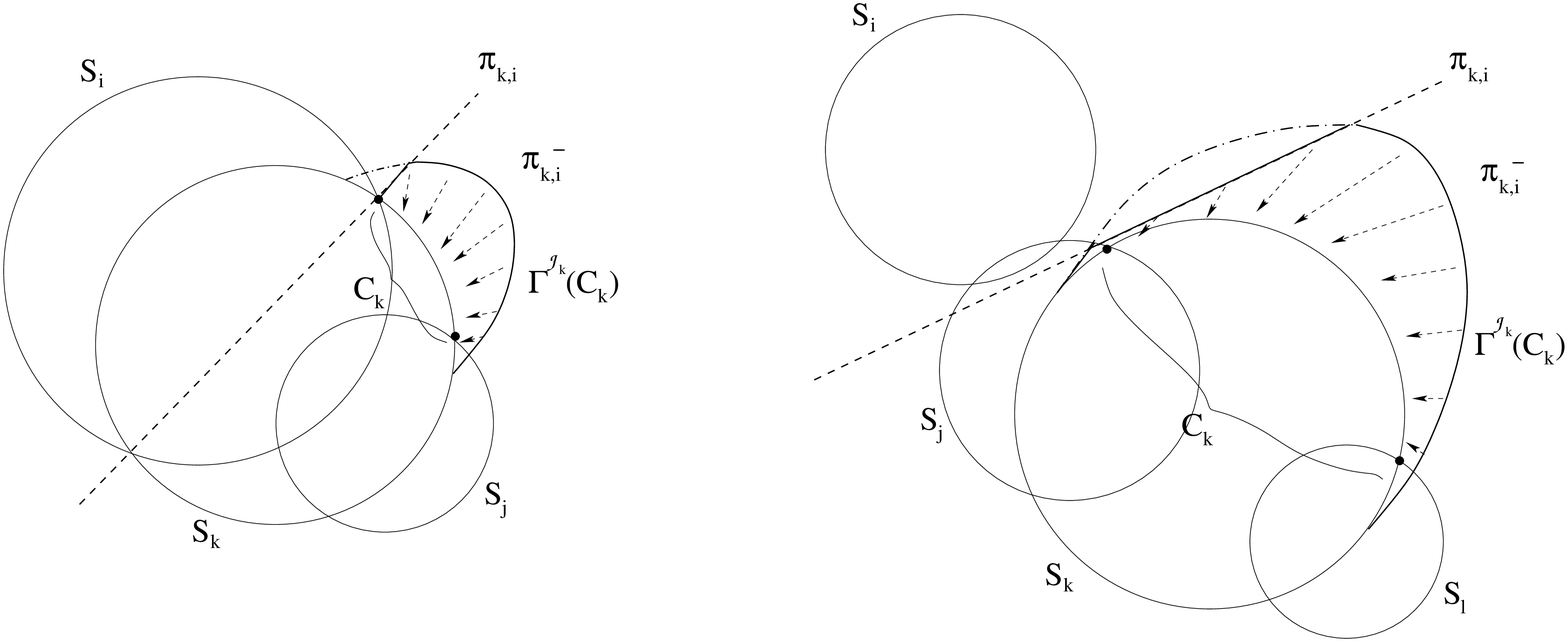,width=15cm}
}}
\vskip 10pt
\centerline{\hbox{\smc figure 3}\hskip 6cm\hbox{\smc figure 4}\hskip 1.3cm}
\endinsert

\n
For every $\J_k\neq\emptyset$, and for any $i\in\J_k$,
call $\pi_{k,i}^-$  the open half space determined by $\pi_{k,i}$ that contains
$C_k\setminus \partial_{\!_{S_k}} {C_k}$.
Then,  setting
%
%
$$
\gather
\Gamma_{_{\overline\ve'}}^{^{\J_k}}
\doteq
\Gamma_{_{\overline\ve'}}^{^{\J_k}}(C_k)
\doteq
\cases
\Gamma_{\overline\ve'}(C_k)\cap
\displaystyle{\bigcap_{i\in\J_k} \pi_{k,i}^-}
\ \ &\text{if}\ \ \ \ \J_k\neq\emptyset\,,
\\
\noalign{\smallskip}
\Gamma_{\overline\ve'}(C_k)
\ \ &\text{otherwise\,,}
\endcases
\qquad\quad k=1,\dots,\nu\,,
\tag 2.38
%
%
%
\\
\noalign{\bigskip}
C\doteq \bigcup_{k=1}^\nu C_k\,,
\qquad\quad
\G \doteq\bigcup_{k=1}^\nu
\overline{\Gamma_{_{\overline\ve'}}^{^{\J_k}}}\,,
\tag 2.39
\\
\noalign{\medskip}
\partial^-\G\doteq\partial\,\G \setminus\bigcup_{k=1}^\nu B_k\,,
\tag 2.40
\endgather
$$
%
one has
$$
\gather
d_C \big(\partial^-\G\big)
\geq c_4\,.
\tag 2.41
\endgather
$$
}
\midinsert
\vskip 2pt
\centerline{\hbox{\psfig{figure=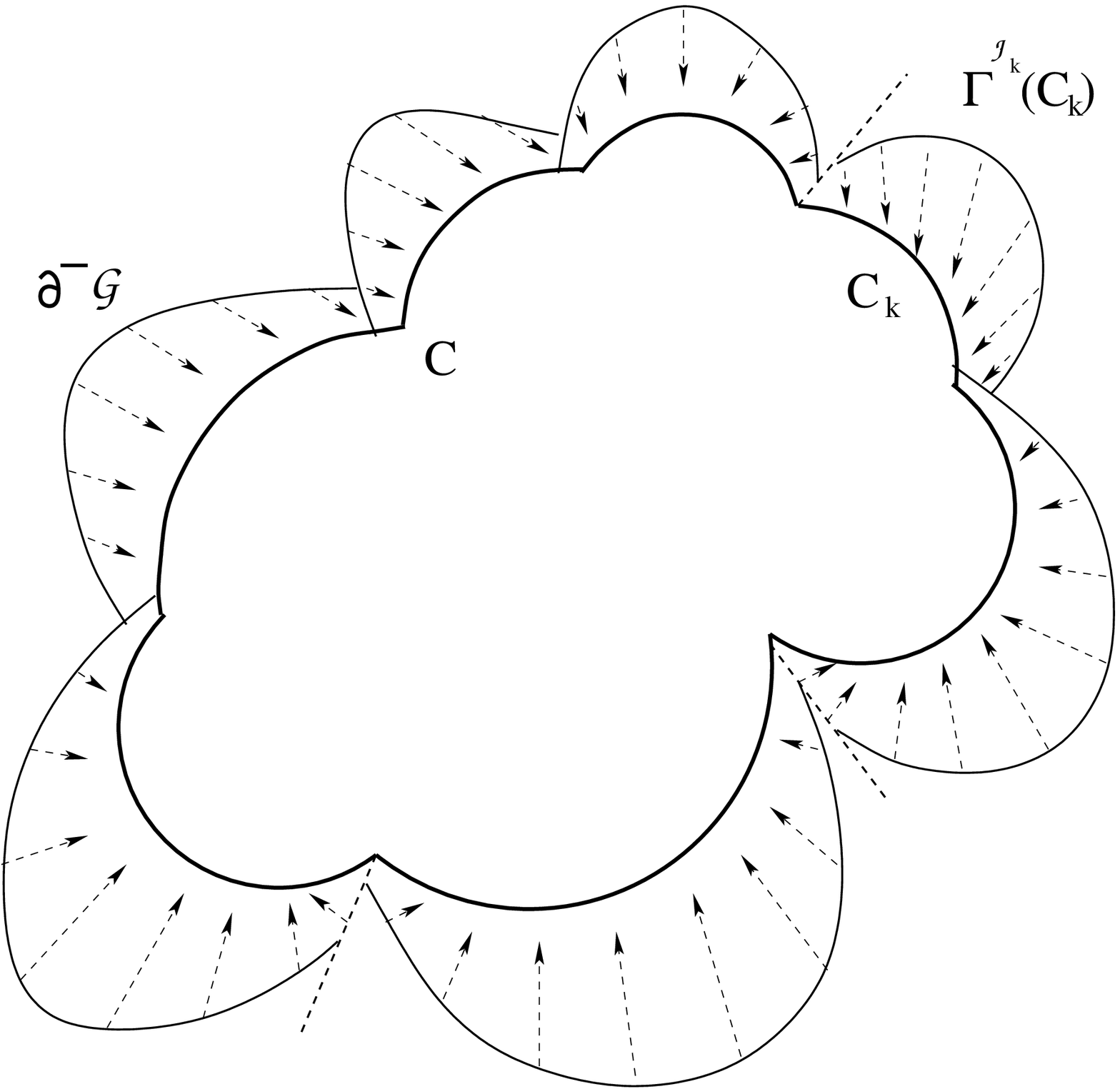,width=8cm}
}}
\vskip 5pt
\centerline{\hbox{\smc figure 5}}
\endinsert

\vs
\n{\medbf 3 - Proof of the theorem}
\v The proof will be given in several steps.
\v

\n{\bf 1.} 
Given $\ve,\,T>0$ ($\ve<\min\{1,T\}$), fix some constant $T'>T+1$, and observe that
by Lemma 1 the value function~$V$ for~(1.18) is Lipschitz continuous
on $\Lambda_{T'}\doteq \big\{ y\in\R^n\,;~~V(y)\leq T'\big\}$.
Hence, by Rademacher's theorem $V$ is
differentiable a.e. in $\Lambda_{T'}$. Then,
letting $\lambda>0$ be the constants provided by Lemma~1
in connection with the set $\Lambda_{T'}$,
for each $y\in\Lambda_{T'}$
at which $V$ is differentiable
define a quadratic function~$V^y$ setting
$$
V^y(x)\doteq V(y)+\nabla V(y)\cdot (x-y)+(1+\lambda)\,|x-y|^2\,.
\tag 3.1
$$
Notice that, because of (2.5), 
there holds 
$$
V(x)+|x-y|^2\leq V^y(x)
\qquad \forall~x\in 
\Lambda_{T'}\,.
\tag 3.2
$$
Moreover, according with Lemma~2, there exists some constant $\ve_0>0$
so that,
given a constant 
$$
0<\ve_1<\frac{\ve}{4 T'}\,,
\tag 3.3
$$
%
for every  $y\in\Lambda_{T'\!,\,\ve_0}\doteq \big\{ y\in\Lambda_{T'}\,;~~|y|\geq\ve_0\big\}$
where $V$ is differentiable
we can choose a control value~$v^y\in\bU$ such that
$$
\big\langle\nabla V(y),\, f(y, v^y)\big\rangle< -1+ \ve_1\,.
\tag 3.4
$$
Choose the constant $\ve_0$ so that, setting
$$
\gathered
\ve'_0\doteq \ve\,\sqrt{\frac{4T'+1}{2(1+2T')}}\,,
\qquad\quad
\ve''_0\doteq \ve\,\sqrt{{2T'\over 1+2T'}}\,,
\\
\noalign{\medskip}
\tau_0\doteq c^{-1}\ln\bigg(\frac{\ve'_0+1}{\ve''_0+1}\bigg),
\endgathered
\tag 3.5
$$
where $c$ denotes the constant in (1.3),
there holds
$$
\ve_0<\min\bigg\{\frac{\ve}{4}\,,\,\ve'_0\bigg\}\,,
\qquad\quad
\frac{\ve_0^2}{\ve^2-\ve_0^2}<\frac{\tau_0}{2}\,.
\tag 3.6
$$
Notice that, by definition (1.17) of $\vp_\ve$,  and because of (3.6),
the value function~$V$ for~(1.18)
satisfies
$$
V(x)\leq \frac{\ve_0^2}{\ve^2-\ve_0^2}<\frac{\tau_0}{2}
\qquad\quad\forall~x\in B_{\ve_0}\,.
\tag 3.7
$$
Next, choose some other constant
$$
L'>L\doteq \text{diam}(\Lambda_{T'}) + 
\sqrt n\,\text{Lip}(V)+\sqrt{n\,(\text{Lip}(V))^2+4T'(1+\lambda)}\,,
\tag 3.8
$$
where $\text{Lip}(V)$ denotes the Lipschitz constant of $V$ on $\Lambda_{T'}$.
Hence, 
since the  assumptions~(H) imply the
Lipschitz continuity in  $x$ of
the function $f(x,u)$ on the compact set~$B_{L'}\times\bU$, uniformly
for~$u\in \bU$, and because also $\nabla V^y$ is Lipschitz continuous
with a Lipschitz constant independent on  $y\in\Lambda_{T'}$,
there will be some constant~$c_5>0$ (depending only on~$L'$)
such that
$$
\aligned
\big|
\big\langle\nabla V^y(x_1),\, f(x_1, u)\big\rangle-
\big\langle\nabla V^y(x_2),\, f(x_2, u)\big\rangle
\big|
&\leq c_5 \cdot\big|x_1-x_2\big|
\\
\noalign{\medskip}
\big|f(x_1, u)-f(x_2, u)\big|
&\leq c_5 \cdot\big|x_1-x_2\big|
\endaligned
\qquad\quad\forall~x_1,x_2\in B_{L'}\,,\ \ u\in \bU\,.
\tag 3.9
$$
Then, setting
$$
c_6\doteq 
\sqrt n\,\text{Lip}(V)+4(1+\lambda)\text{diam}(\Lambda_{T'}) \,,
\tag 3.10
$$
and choosing $\ve_2>0$ so that
$$
\ve_2<\min\bigg\{\sqrt{\frac{\ve}{3}},\,\frac{\ve_1}{8c_5c_6},
\,\sqrt{\frac{\tau_0}{2}},\,\frac{T'-T}{1+\text{Lip}(V)},\,L'-L\bigg\}\,,
\tag 3.11
$$
we deduce from (3.4), (3.9) that, 
for every $y\in\Lambda_{T'\!,\,\ve_0}$
where $V$ is differentiable
there holds
$$
\big\langle\nabla V^y(x),\, f(x, v^y)\big\rangle < -1+ 2\ve_1
\qquad\quad
\forall~x\in B(y,2\ve_2)\cap B_{L'}\,.
\tag 3.12
$$
\v
\n{\bf 2.}
By the Lipschitz continuity of $V$ on the set $\Lambda_{T'}$ it follows that,
for each $y\in\Lambda_{T'\!,\,\ve_0}$
at which $V$ is differentiable,
there holds
$$
\big|V^y(x)-V(x)\big|\leq c_7\cdot \big|x-y\big|
\qquad\quad\forall~x\in\Lambda_{T'}\,,
$$
for some positive constant $c_7$\,.
Hence, since the set $\Lambda_{T',\ve_0}$ is compact
(cfr. point~{\bf 1} of the proof of Lemma~1), we can cover 
it with finitely many balls (of sufficiently small radius),
centered at points of $\Lambda_{T'\!,\,\ve_0}$ where $V$ is differentiable, 
say $y_1,\ldots, y_N$,
so that,  setting
$$
\aligned
V_i(x)&\doteq V^{y_i}(x)\,,
\quad
1\leq i \leq N\,,%
\\
\noalign{\medskip}
\tV(x)&\doteq \min_i V_i(x)
\endaligned
\qquad\quad\forall~x\in\R^n\,,
\tag 3.13
$$
there holds
$$
V(x)\leq\tV(x)\leq V(x)+\ve_2^2
\qquad\quad\forall~x\in\Lambda_{T'\!,\,\ve_0}\,.
\tag 3.14
$$
Next, observing that (3.2) implies
$$
V_i(x)>V(x)+\ve_2^2
\qquad\quad\forall~x\in\Lambda_{T'}\setminus B(y_i,\ve_2)\,,
$$
we deduce from (3.14)
that
$$
\tV(x)<V_i(x)
\qquad\quad\forall~x\in\Lambda_{T'\!,\,\ve_0}\setminus B(y_i,\ve_2)\,.
\tag 3.15
$$
Relying on (3.15), letting
$$
\P_i\doteq \big\{  x\in\R^n\,;~~ V_i(x)= \tV(x)\big\}\,,
\tag 3.16
$$
we find that
$$
\P_i\cap
\Lambda_{T'\!,\,\ve_0}\subset
B(y_i,\ve_2)\,,
\qquad\quad 1\leq i \leq N\,.
\tag 3.17
$$
Hence, by (3.12), (3.17) we have
$$
\big\langle\nabla V_i(x),\, f(x, v^i)\big\rangle < -1+ 2\ve_1
\qquad\quad\forall~x\in B\big(\P_i,\,\ve_2\big)\cap \Lambda_{T'\!,\,\ve_0}
\cap B_{L'},\ \ 1\leq i \leq N\,,
\tag 3.18
$$
where we have set $v^i\doteq v^{y_i}$ ($1\leq i \leq N$),
while (3.2) yields
$$
\tV(x)\geq V(x)+|x-y_i|^2
\qquad\quad\forall~x\in \P_i\cap\Lambda_{T'},\ \ 1\leq i \leq N\,.
\tag 3.19
$$
\v

\n{\bf 3.} The patchy feedback $u=U(x)$ will be constructed
looking at the level sets of the function $\tV$
defined in (3.13).
To this end, observe first that, because of (3.14),
and by the choice~(3.11) of~$\ve_2$, there holds
$$
\tV(x)<T'-\text{Lip}(V)\cdot\ve_2\qquad\quad\forall~x\in\Lambda_{T}\,.
\tag 3.20
$$
Moreover, notice that, relying on the definitions (3.1) of $V^{y_i}$
and (3.8) of the constant $L$, one finds
$$
\aligned
1\leq i\leq N\,,\ \ V_i(x)<T'\qquad
\Longrightarrow
\qquad\
|x|&\leq |x-y_i|+\text{diam}(\Lambda_{T'})
\\
\noalign{\medskip}
&\leq
|\nabla V(y_i)|+
\sqrt{|\nabla V(y_i)|^2+4T'(1+\lambda)}+\text{diam}(\Lambda_{T'})
\\
\noalign{\medskip}
&\leq
\sqrt n\,\text{Lip}(V)+\sqrt{n\,(\text{Lip}(V))^2+4T'(1+\lambda)}+\text{diam}(\Lambda_{T'})
\\
\noalign{\medskip}
&\leq L\,,
\endaligned
\tag 3.21
$$
which, in turn, by the definition (3.13) of $\tV$, and because of (3.14),  yields
$$
\big\{
x\in\R^n~;~~\tV(x)<T'
\big\}
\subset \Lambda_{T'}\cap B_L\,.
\tag 3.22
$$
On the other hand, observe that all level sets of each quadratic function~$V_i$ are spheres.
Therefore, 
every level set
$$
\Sigma_\tau\doteq \big\{ x\in\R^n~;~~\tV(x)=\tau\big\}
\qquad\tau\geq \tau_0\,,
\tag 3.23
$$
is contained in a finite union of spheres, and  each upper level set
$\big\{ x\in\R^n;~~\tV(x)\geq\tau\big\}$, $\tau\geq\tau_0$, is connected.
Moreover, notice that
by (3.7), (3.11), (3.14), we derive
$$
\max_{|x|=\ve_0}\tV(x)<\frac{\tau_0}{2}+\ve_2^2\leq \tau_0\,,
\tag 3.24
$$
and hence we find that
$$
\big\{ x\in\R^n;~~\tV(x)\geq\tau_0\big\}\cap B_{\ve_0}=\emptyset\,.
\tag 3.25
$$
Thus,  setting
%
$$
T''\doteq T'-\text{Lip}(V)\cdot\ve_2\,,
\qquad\quad
\D\doteq\big\{
x\in\R^n~;~~\tau_0<\tV(x)<T''
\big\}\,,
\tag 3.26
$$
thanks to  (3.11), (3.14), (3.20), (3.22), (3.25),
we deduce that
$$
\Lambda_{T,\ve}\subset B\big(\D,\,\ve_2\big)
\subset\Lambda_{T'\!,\,\ve_0}\cap B_L\,.
\tag 3.27
$$
We will establish the theorem by constructing the
patchy feedback $u=U(x)$ on the domain $\D$.
Notice that, with the same arguments used in the proof of Lemma~2,
by the choice of the constants~$\ve'_0, \, \tau_0$ in (3.5) we find that
$$
V(x)>\tau_0\qquad\quad\forall~x\in\Lambda_{T'}\,, \quad |x|\geq \ve'_0\,.
\tag 3.28
$$
Hence, since the definition (3.13) of $\tV$  implies 
$$
V(x)\leq \tau_0\qquad\quad\forall~x\in\Sigma_{\tau_0}\,,
$$
we deduce from (3.6), (3.28)  that 
$$
\Sigma_{\tau_0}\subset B_{\ve'_0}\subset B_{\ve}\,.
\tag 3.29
$$
Next, observe that, since all functions $V_i$, $1\leq i \leq N$, have the same coefficient
of the quadratic term,
it follows that,
for each couple of indices $k\not= i$, the set
$$
\pi_{k,i}\doteq\pi_{i,k}\doteq \big\{ x\in\R^n\,;~~V_k(x)=V_i(x)\big\}
\tag 3.30
$$
is an hyperplane, and the difference of the gradients
$\nabla V_i(x)-\nabla V_k(x)$ is a constant vector on $\pi_{k,i}$.
Then, letting
$\bfn_{k,i}$ denote the unit normal to $\pi_{k,i}$,
pointing towards the half space 
$$
\pi_{k,i}^+\doteq\pi_{i,k}^-\doteq
\big\{ x\in\R^n\,;~~V_k(x)>V_i(x)\big\}\,,
 \tag 3.31
$$
%
%
one has
$$
\nabla V_i(x)-\nabla V_k(x)= -c\,\bn_{k,i}
\qquad\quad\forall~x\in\pi_{k,i}\,,
\tag 3.32
$$
for some constant $c=c_{k,i}\geq 0$.
Denote as $\pi_{k,i}^-$ the other half space 
determined by $\pi_{k,i}$, i.e. set
$$
\pi_{k,i}^-\doteq\pi_{i,k}^+\doteq
\big\{ x\in\R^n\,;~~V_k(x)<V_i(x)\big\}\,.
 \tag 3.33
$$

\v
\n{\bf 4.} 
The basic step in the construction of $U(x)$ is the following. 
We shall fix a suitably small time size~$\Delta t$ and, in connection with 
an increasing sequence of times $\{\tau_m\}_{m\geq 0}$ with 
the property 
$$
\exists~p\quad\text{s.t.}\quad\tau_{m+p}>\tau_m+\Delta_t
\qquad\forall~m\,,
$$
we will construct,
for every $m\geq 0$, a patchy feedback whose domain
contains the region
$$
\D_m\doteq\big\{x\in\R^n\,;~~\tau_m< \tV(x)\leq\tau_{m+1}\big\}\,,
$$
so that all the 
trajectories $x(t)$ of the corresponding closed-loop system (1.2)
satisfy
$$
\frac{d}{dt} \tV\big(x(t)\big)
\leq -1+3 \ve_1\qquad\ \text{for \ \ \ a.e.}~t\,,
$$
and eventually enter the set where $\tV<\tau_m$.
To this end, fix any $\tau\in[\tau_0, T'[$ and consider the level set~$\Sigma_\tau$ of~$\tV$.
By construction, $\Sigma_\tau$
is contained in the union of finitely
many spheres, say 
$S_{i_1},\ldots,S_{i_{\nu_\tau}}$. Here we denote as $S_{i_\ell}\doteq\big\{ x\in\R^n\,;~~
V_{i_\ell}(x)=\tau\big\}$ the surface of the ball
$B_{i_\ell}\doteq\big\{ x\in\R^n\,;~~
V_{i_\ell}(x)\leq\tau\big\}$.
Notice that, since the definition (3.13) of $\tV$
implies $\tV(x)<\tau$ for all
$x\in\overset\circ\to B_{i_\ell}$,
by definition (3.23) it follows that
$$
\gather
\Sigma_\tau=\bigcup_{\ell=1}^{\nu_\tau} \Sigma_{\tau,i_\ell}\,,
\qquad\quad
\Sigma_{\tau,i_\ell}\doteq S_{i_\ell}\setminus
\bigcup_{q=1}^{\nu_\tau} \overset\circ\to B_{i_q}\,,
\tag 3.34
\\
\noalign{\medskip}
\big\{ x\in\R^n\,;~~
\tV(x)<\tau\big\}=
\bigcup_{\ell=1}^{\nu_\tau}\overset\circ\to B_{i_\ell}\,.
\tag 3.35
\endgather
$$
We can assume that the set of indices
$\I_\tau\doteq \{ i_1,\ldots, i_{\nu_\tau}\}$
includes 
only those indices $i\in\{1,\ldots,N\}$ 
for which there exists some point $\overline x\in\D$
satisfying
$$
\tau=V_i(\overline x)<V_j(\overline x)
\qquad\forall~ j\neq i.
$$
This means that
$$
\I_\tau=\Big\{
i\in\{1,\dots,N\}\,;\ 
\Big(\P_i \setminus \bigcup_{j\neq i} \P_j\Big)\cap \Sigma_\tau
\neq \emptyset
\,\Big\}\,,
\tag 3.36
$$
and, in particular, implies that
$$
S_{i}\setminus
\bigcup_{j\in\I_\tau} \overset\circ\to B_{j}
\neq \emptyset\qquad\quad\forall~i\in\I_\tau\,.
\tag 3.37
$$
Moreover, we may write $\Sigma_\tau$ as the union of $\eta_\tau$
connected components $\Sigma_\tau^{\strut 1},\,\dots,\, \Sigma_\tau^{\strut \eta_\tau}$, so that
setting
$$
\I_\tau^h\doteq
\Big\{
i\in\I_\tau\,;\ 
\Big(\P_i \setminus \bigcup_{j\neq i} \P_j\Big)\cap \Sigma_\tau^{\strut h}
\neq \emptyset
\,\Big\}\,,
\tag 3.38
$$
there holds
$$
S_{i}\cap
\bigcup_{j\in\I^h_\tau} \overset\circ\to B_{j}
\neq \emptyset\qquad\qquad\forall~i\in\I_\tau^h\,,\qquad h=1,\dots,\eta_\tau\,.
\tag 3.39
$$
Notice also that, by (3.13), (3.23), (3.34), (3.37), 
every set $\Sigma_{\tau,i}$, $i\in\I_\tau$ is nonempty and one has
$$
\gathered
\Sigma_{\tau,i}=\big\{ x\in\R^n\,;~~
V_i(x)=\tV(x)=\tau\big\}
\endgathered
\qquad\quad
\forall~i\in\I_\tau\,,
\tag 3.40
$$
while the definitions  (3.30), (3.31), (3.33) imply
$$
\gathered
\pi_{k,i}\cap S_k=\pi_{k,i}\cap S_i=S_k\cap S_i\,,
\\
\noalign{\medskip}
\Sigma_{\tau,k}\subset \pi_{k,i}\cup\pi_{k,i}^-\,,
\qquad\qquad
\Sigma_{\tau,i}\subset \pi_{k,i}\cup\pi_{k,i}^+
\endgathered
\qquad\quad
\forall~k,i\in\I_\tau\,.
\tag 3.41
$$
Therefore, relying on (3.41) we deduce that, for every pair of indices
$k,i\in\I_\tau$, $k\neq i$, one of the following two cases occurs:
\item{--} if $\overset\circ\to B_k\cap \overset\circ\to B_i\neq\emptyset$\,,
then $\pi_{k,i}$ is the hyperplane passing through $S_k\cap S_i$;

\item{--}  if $\overset\circ\to B_k\cap \overset\circ\to B_i=\emptyset$\,,
then $S_k\subset \pi_{k,i}\cup\pi_{k,i}^-$ and  $S_i\subset \pi_{k,i}\cup\pi_{k,i}^+$, i.e.
$\pi_{k,i}$ is an hyperplane separating~$S_k$ and~$S_i$.

\v

\n{\bf 5.} By the above construction, 
and relying on (3.11), (3.22), we deduce
that
$$
\gathered
\Sigma_{\tau,i}
\subset \P_i\cap B_L\,,
\\
\noalign{\medskip}
B\big(\Sigma_{\tau,i},\,\ve_2\big)
\subset B_{L'}\,,
\endgathered
\qquad\ \forall~\tau\in[\tau_0,T''[\,,\ \ i\in\I_\tau\,.
\tag 3.42
$$
Hence, thanks to (3.18), (3.27), (3.42),
we find
$$
\big\langle\nabla V_i(x),\, f(x, v^i)\big\rangle< -1+ 2\ve_1
\qquad\ \forall~x\in B\big(\Sigma_{\tau,i},\ve_2\big)\,,
\quad \tau\in[\tau_0,T''[\,,\ \ i\in\I_\tau\,.
\tag 3.43
$$
Relying on (3.43), 
we shall construct around each set $\Sigma_{\tau,i}$, $i\in\I_\tau$,
a lens-shaped domain~$\Gamma_{\tau,i}$
of the form~(2.32) as in Lemma 3, 
so that the boundary of 
 $\Gamma_{\tau,i}$ is transversal to the flow
of the vector field~$g_i(x)\doteq f(x, v^i)$.
Namely, letting $x(t;y, v^i)$ denote the solution of the
Cauchy problem $\dot x =g_i(x)$, $x(0)=y$, we will prove the following
\v

\n
{\bf Claim~1.} {\it There exists a positive constants $\ve_3$
so that, for every given $\tau\in[\tau_0,T''[\,,$
$k\in\I_\tau$, the vector field $g_{\tau,k}(x)\doteq f(x, v^k)$
is tranversal to 
the  boundary of the domain
$$
\Gamma_{\tau,k}\doteq
\Big\{ x(s;y, v^k)\,; ~~y\in B\big(\Sigma_{\tau,k},\, \ve_3\big)\cap S_k\,,~~
{d_{_{\Sigma_{\tau,k}}}}^{\!\!\!\!\!\!\!\!\!2} \ \ (y)-\ve_3^2< s\leq 0\,
\Big\}\,.
\tag 3.44
$$
Namely, it points strictly inward on the upper boundary
$$
\partial^-\,\Gamma_{\tau,k}\doteq \partial\,\Gamma_{\tau,k}
\setminus 
\bigcup_{j\in\I_\tau} B_j
$$
and strictly outward on 
the lower boundary
$$
\partial^+\,\Gamma_{\tau,k}\doteq \partial\,\Gamma_{\tau,k}
\cap S_k\,.
$$
Moreover,  there holds
$$
|g_{\tau,k}(x)| \geq c_8\qquad\quad\forall~x\in\overline{\Gamma_{\tau,k}}\,,
\tag 3.45
$$

$$
\overline{\Gamma_{\tau,k}}\subset B\big(\Sigma_{\tau,k},\,\ve_2\big)\subset B\big(y_k,\,2\ve_2\big)\,,
\tag 3.46
$$
for some constant $c_8>0$ independent on $\tau\in[\tau_0,T''[\,,$
$k\in\I_\tau$.
}
\v

\n{\bf 6.}\ {\smc Proof of Claim~1.} 
In order to establish the claim, 
we shall first derive an upper and lower 
 uniform bound for the radii of 
the spheres 
$$
S_i\doteq \big\{ x\in\R^n\,;~~ V_i(x)=\tau\big\}\qquad\
i\in\I_\tau\,, \quad \tau\in[\tau_0,T''[\,,
\tag 3.47
$$
and we will prove an
a priori estimate for $\langle\bfn_{i},\, f(x, v^i)\rangle$,
$x\in\Sigma_{\tau,i}$,
 independent of $\tau\in[\tau_0,T''[$\,, \linebreak and $i\in\I_\tau$ 
($\bfn_{i}$ denoting the unit outer normal to $S_i$).
%
To this end observe that by (1.3)
one has 
$$
\big|g_{\tau,i}(x)\big|=
\big|f(x, v^i)\big|\leq c_9\doteq c\big(1+L')
\qquad\forall~x\in B_{L'}\,.
\tag 3.48
$$
Then, for every fixed $i\in\I_\tau$, $\tau\in[\tau_0,T''[$\,,
writing
$$
V_i(x)=(1+\lambda)\,|x-\omega_i|^2 + b_i
\qquad\quad\forall~x\in S_i\,,
$$
for some point $\omega_i\in\R^n$ and some constant $b_i$,
and using  (3.3), (3.42), (3.43), (3.48), we derive the estimate
$$
\aligned
2(1+\lambda)|x-\omega_i|= \big|\nabla V_i(x)\big|
&\geq \frac{\big|\big\langle\nabla V_i(x),\, f(x, v^i)\big\rangle\big|}{\big|f(x, v^i)\big|}
\\
\noalign{\medskip}
&\geq{1\over 2 c_9}
\qquad\ \ \forall~x\in\Sigma_{\tau,i}\,.
\endaligned
\tag 3.49
$$
On the other hand,
from the definition (3.1) of $V_i=V^{y_i}$,
recalling that $y_i\in\Lambda_{T'}$,
and relying on~(3.10), (3.27), 
one deduces the a-priori bound
$$
\big|\nabla V_i(x)\big|\leq 
\big|\nabla V(x)\big|+4(1+\lambda)\text{diam}(\Lambda_{T'}) 
\leq c_6
\qquad\ \ \forall~x\in\D\,.
\tag 3.50
$$
Hence, thanks to (3.49), (3.50),
we find that the radius $r_i=|x-\omega_i|$, $x\in \Sigma_{\tau,i}$\,, of the sphere $S_i$ satisfies
$$
{1\over 4c_9(1+\lambda)}
\leq r_i\leq
{{c_6}\over {2(1+\lambda)}}
\qquad\quad
\forall~i\in\I_\tau\,, \quad \tau\in[\tau_0,T''[\,,
\tag 3.51
$$
while (3.3), (3.43), together with (3.50), yield
$$
\big\langle\bfn_i ,\, f(x, v^i)\big\rangle=
{\big\langle\nabla V_i(x),\, f(x, v^i)\big\rangle\over 
\big|\nabla V_i(x)\big|}\leq -{1\over 2c_6}
\qquad\ \ \forall~x\in\Sigma_{\tau,i}\,,
\quad i\in\I_\tau\,, \quad \tau\in[\tau_0,T''[\,.
\tag 3.52
$$
Therefore, because of (3.51), (3.52),
we can apply Lemma~3 to every set $\Sigma_{\tau,k}$,
$\tau\in[\tau_0,T''[$\,, $k\in\I_\tau$,
in connection with the vector field $g_{\tau,k}$.
Thus we deduce the existence of some constant~$\ve_3>0$,  
so that the field~$g_{\tau,k}(x)=f(x,v^k)$ is transversal to the boundary of
the domain $\Gamma_{\tau,k}$ defined in~(3.44).
Concerning (3.46), observe that choosing $\ve_3$
such that $\ve_3 (c_9\,\ve_3+1)<\ve_2$, thanks to~(3.42), (3.48)
we obtain
$$
\aligned
d_{\Sigma_{\tau,k}}\big(
x(s;\,y, v^k)\big)+d_{\Sigma_{\tau,k}}(y)&\leq
\big|x(s;\,y, v^k)-y\big|
\\
\noalign{\medskip}
&\leq |s|\,\|g_{\tau,k}\|_{\strut\L^\infty(B_{L'})}+d_{\Sigma_{\tau,k}}(y)
\\
\noalign{\smallskip}
&\leq |s|\,c_9+\ve_3
< \ve_2
\qquad\quad \forall~y\in B\big(\Sigma_{\tau,k},\, \ve_3\big)\cap S_k\,,\quad
-\ve_3^2< s\leq 0\,.
\endaligned
\tag 3.53
$$
Hence,  because of (3.42), (3.17), (3.27), relying on (3.53)
we find
$$
\overline{\Gamma_{\tau,k}}\subset B\big(\Sigma_{\tau,k},\ve_2\big)
\subset B\big(\P_k,\ve_2\big)\subset B\big(y_k,2\ve_2\big)\,,
\tag 3.54
$$
which
proves (3.46).
Finally, observe that,
for every given $k\in\I_\tau\,, \tau\in[\tau_0,T''[$\,,
fixing some \linebreak point~$\overline x\in \Sigma_{\tau,k}$,
thanks to (3.46), 
and because of (3.3), (3.9), (3.10), (3.11), (3.52), we derive
$$
\aligned
\big|f(x, v^k)\big|&\geq
\big|f(\overline x, v^k)\big|-
\big|f(x, v^k)-f(\overline x, v^k)\big|
\\
\noalign{\medskip}
&\geq \big|\big\langle\bfn_k ,\, f(\overline x, v^k)\big\rangle \big|-
c_5 \cdot\big|x-\overline x\big|
\\
&\geq{1\over 2c_6}-4c_5\cdot\ve_2
\\
&\geq{1\over 4c_6}\qquad\ \ 
\forall~x\in\overline{\Gamma_{\tau,k}}\,,
\endaligned
\tag 3.55
$$
which yields (3.45), thus completing the proof of our claim.

\v

\n{\bf 7.} Given $\tau\in[\tau_0, T''[$\,, $k\in\I_\tau$,
consider now the domain $\Gamma_{\tau,k}$
defined in (3.44), and observe that, because of~(3.16), (3.43), (3.46),
every trajectory $x(t)$ of $\dot x = g_{\tau,k}(x)$,
passing through points \linebreak of~$\Gamma_{\tau,k}\cap \P_k$,
satisfies 
$$
\aligned
\tV\big(x(t)\big)&= V_k\big(x(t)\big)
\\
\noalign{\smallskip}
&=V_k\big(x(s)\big)+\int_s^t
\big\langle\nabla V_k(x(\sigma)),\, f(x(\sigma), v^k)\big\rangle~d\sigma
\\
\noalign{\medskip}
&\leq V_k\big(x(s)\big)+(-1+2\ve_1)(t-s)
\\
\noalign{\medskip}
&=
\tV\big(x(s)\big)+(-1+2\ve_1)(t-s)
\qquad\ \ \ \forall~
t>s\,.
\endaligned
\tag 3.56
$$
However,  there may well be points $x(t)\in \Gamma_{\tau,k}$ where
$V_k(x(t))>\tV(x(t))$. Near these points there is no guarantee that
(3.56) should hold.   To address this difficulty, we will
consider the set
of all indices $i\not= k$ 
such that $V_i(\overline x)<V_k(\overline x)$ for some $\overline x\in \Gamma_{\tau,k}$\,,
and such that
$$
\min
\Sb
\\
\noalign{\vskip 1pt}
x\in \overline{\Gamma_{\tau,k}}
\endSb
\big\langle\nabla V_k(x)-\nabla V_i(x),\, f(x, v^k)\big\rangle<0\,.
\tag 3.57
$$ 
In this case, we shall replace $\Gamma_{\tau,k}$ with the smaller domain
$$
\Gamma_{\tau,k}\bigcap \big\{ x\in\R^n\,;~~V_k(x)<V_i(x)\big\}\,.
$$
Then, setting
$$
\J_{\tau,k}
\doteq
\Big\{
i\in\{1,\dots,N\}\setminus\{k\}\,;\ 
\P_i\cap\Gamma_{\tau,k}\neq\emptyset\,, \ \
\min
\Sb
\\
\noalign{\vskip 1pt}
x\in 
\overline{\Gamma_{\tau,k}}
\endSb\big\langle\nabla V_k(x)-\nabla V_i(x),\, f(x, v^k)\big\rangle<0
\,\Big\}\,,
\tag 3.58
$$
%
%
consider the domain
$$
\widetilde\Gamma_{\tau,k}\doteq
\cases
\Gamma_{\tau,k}\cap
\displaystyle{\bigcap_{i\in\J_{\tau,k}} \pi_{k,i}^-}\quad &\text{if}\qquad \J_{\tau,k}\neq\emptyset\,,
\\
\noalign{\smallskip}
\Gamma_{\tau,k}
\quad &\text{otherwise,}
\endcases
\tag 3.59
$$
which, according with the definitions (2.32), (2.38), is precisely
equal to~$\Gamma_{{\ve_3}}^{\!^{\J_{\tau,k}}}(\Sigma_{\tau,k})$.
\midinsert
\vskip 10pt
\centerline{\hbox{\psfig{figure=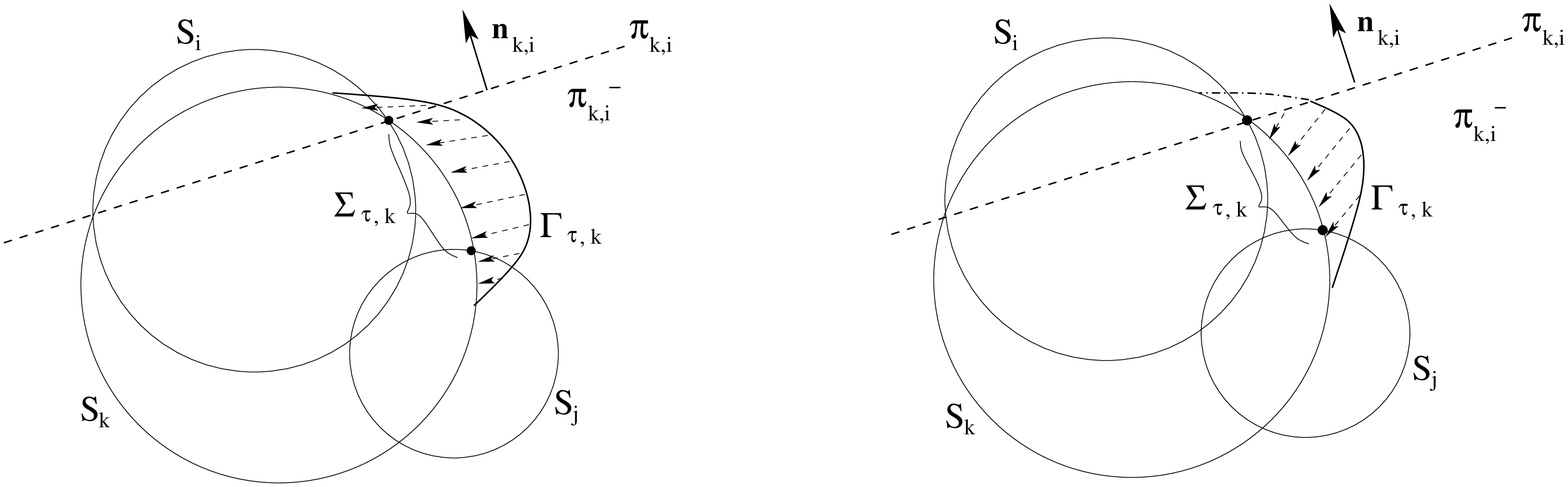,width=15cm}}}
\vskip 10pt
\centerline{\hbox{\smc figure 6} \ \ ($i\notin\J_{\tau,k}$)
\hskip 5cm\hbox{\smc figure 7} \ \ ($i\in\J_{\tau,k}$)}
\vskip 5pt
\endinsert
\n
Notice that, because of (3.37), (3.39), (3.51), (3.52),
and by the observations at point {\bf 4}, for every \linebreak fixed
$h=1,\dots, \eta_\tau$, the spheres $S_i$, $i\in\I^h_\tau$,  
and the collection of hyperplanes
$$
\Pi^h_\tau\doteq\big\{\pi_{k,i}\,;~~k,i\in\I^h_\tau\big\}
$$
%
satisfy the assumptions
of Lemma~4. 
Hence, in the case where
$$
\bigcup_{k\in\I^h_\tau}\J_{\tau,k}\subset\I^h_\tau\,,
$$
we are in the position to apply the conclusion of Lemma~4
 in connection with the collection of hyperplanes $\Pi^h_\tau$
and of sets
$$
\big\{\Sigma_{\tau,i}\,;~~i\in\I^h_\tau\big\}
$$
defined in (3.34),
in order to derive a uniform  estimate of the distance of the (upper) boundary of 
$$
\G^h_\tau\doteq \bigcup_{k\in\I^h_\tau}\overline{\widetilde\Gamma_{\tau,k}}
\tag 3.60
$$
from the set
$$
\Sigma^h_\tau=\displaystyle{\bigcup_{i\in\I_\tau^h}\Sigma_{\tau,i}}\,.
$$
As a consequence, we obtain an estimate
of the decrease of $\tV$ along trajectories of~$g_{\tau,k}$ passing
through $\widetilde\Gamma_{\tau,k}$, $k\in\I^h_\tau$. More precisely, setting
$$
\gathered
\I_1^\ast\doteq
\Big\{\tau\in[\tau_0, T''[~;~~\bigcup_{k\in\I^h_\tau}\J_{\tau,k}\subset\I^h_\tau\qquad
\forall~h=1,\dots, \eta_\tau
 \Big\}\,,
\\
\noalign{\medskip}
\G_\tau\doteq \bigcup_{k\in\I_\tau}\overline{\widetilde\Gamma_{\tau,k}}
=\bigcup_{h=1}^{\eta_\tau} \G^h_\tau\,,
\endgathered
\tag 3.61
$$
we will prove the following
\midinsert
\vskip 10pt
\centerline{\hbox{\psfig{figure=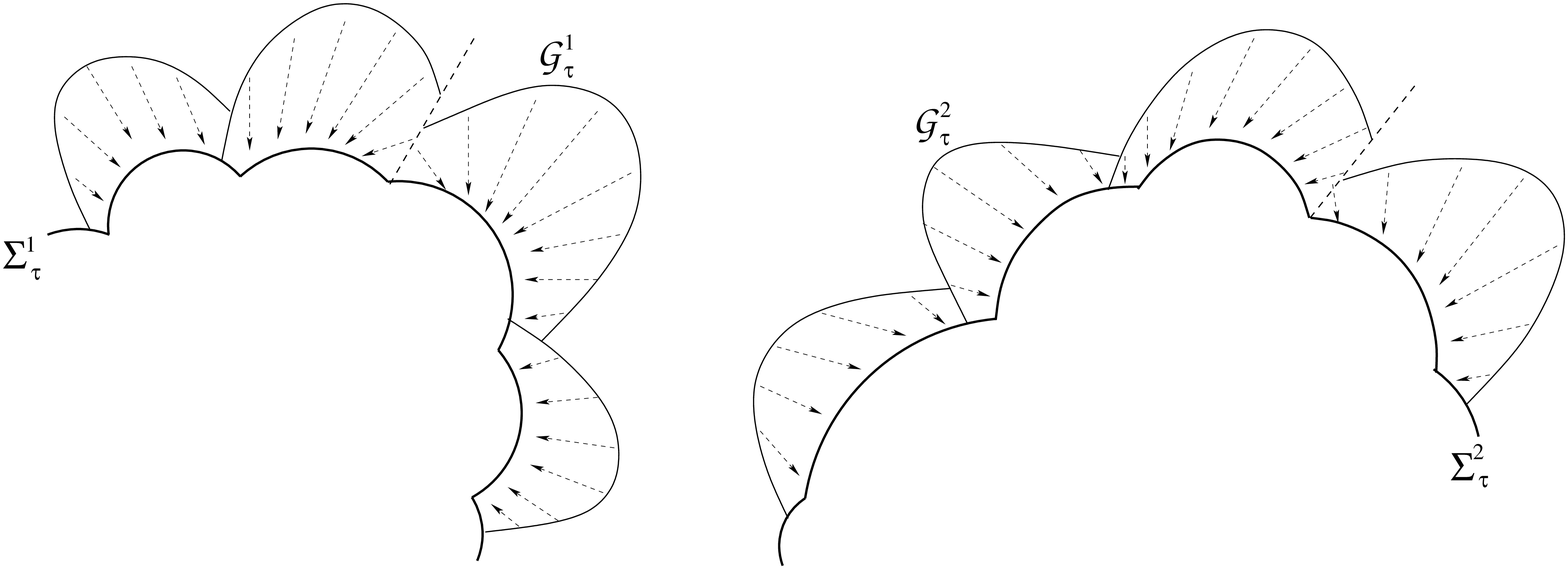,width=14cm}}}
\vskip 10pt
\centerline{\hbox{\smc figure 8}}
\vskip 5pt
\endinsert
\v

\n
{\bf Claim~2.} {\it 
The domains $\widetilde\Gamma_{\tau,k}$, $k\in\I_\tau$, $\tau\in[\tau_0,T''[\,,$ defined in (3.59)
enjoy the following properties.
\item{i)} For any $k\in\I_\tau$, $\tau\in[\tau_0,T''[\,,$ the 
vector field $g_{\tau,k}(x)\doteq f(x, v^k)$  points strictly inward
at every point of the upper boundary 
$$
\partial^-\widetilde\Gamma_{\tau,k}\doteq
\partial\widetilde\Gamma_{\tau,k}\setminus \bigcup_{j\in\I_\tau} B_j\,.
\tag 3.62
$$
\item{ii)}
For any $y\in \displaystyle{\overline{\widetilde\Gamma_{\tau,k}}\setminus\bigcup_{j\in\I_\tau} B_j}$, 
$k\in\I_\tau$, $\tau\in[\tau_0,T''[\,,$ there exists a time $\T_{\tau,k}(y)>0$
so that one has
%
$$
x\big(\T_{\tau,k}(y);\,y,v^k\big)\in\Sigma_\tau\,,
\tag 3.63
$$
%
%
$$
x(t;\,y,v^k)\in \widetilde\Gamma_{\tau,k}
\qquad\forall~t\in\,]0,\,\T_{\tau,k}(y)]\,,
\tag 3.64
$$
and there holds
$$
\tV\big(x(t;\,y,v^k)\big)\leq
\tV\big(x(s;\,y,v^k))+(-1+3\ve_1)(t-s)
\qquad\quad\forall~0\leq s<t\leq \T_{\tau,k}(y)\,,
\tag 3.65
$$
where $\ve_1$ is the constant satisfying (3.3).
\medskip
\item{iii)} 
%
For any $\tau\in[\tau_0, T''[$\,, one has
$$
\tau^\sharp\doteq
\sup \Big\{t\in[\tau,\,T''[~;\quad \Sigma_s\subset \displaystyle{\overset\circ\to \G_\tau\setminus \bigcup_{j\in\I_\tau} 
B_j}\quad\forall~s\in\,]\tau,\,t]\Big\}>\tau\,.
\tag 3.66
$$
Moreover, there exists a positive constant $\ve_4$ so that
there holds
$$
\tau^\sharp>\tau+\ve_4
\qquad\quad\forall~\tau\in\I_1^\ast\,.
\tag 3.67
$$
}

\v

\n{\bf 8.}\ {\smc Proof of Claim~2.} 
By Claim~1 we know that, for every
$k\in\I_\tau$, the vector field $g_{\tau,k}$ is inward-pointing on the 
region~$\partial^-\,\widetilde\Gamma_{\tau,k}\cap\partial^-\,\Gamma_{\tau,k}$.
On the other hand, recalling (3.32), the inequality~(3.57)
guarantees that $g_{\tau,k}$ enjoys the inward-pointing condition also
at the boundary points~$x\in\partial^-\,\widetilde\Gamma_{\tau,k}
\cap\overset\circ\to\Gamma_{\tau,k}\cap\pi_{k,i}$, $i\in\J_{\tau,k}$.
Then, observing that 
$$
\partial^-\,\widetilde\Gamma_{\tau,k}\setminus \partial^-\,\Gamma_{\tau,k}=
\partial^-\,\widetilde\Gamma_{\tau,k}\cap \overset\circ\to\Gamma_{\tau,k}
\cap \displaystyle{\bigcup_{i\in\J_{\tau,k}}\pi_{k,i}}\,,
$$
by continuity it follows that $g_{\tau,k}(x)\in\overset\circ\to T_{\widetilde\Gamma_{\tau,k}}(x)$
at every  point $x\in\partial^-\,\widetilde\Gamma_{\tau,k}$ 
($\,\overset\circ\to T_{\widetilde\Gamma_{\tau,k}}$
denoting the interior of the tangent cone to $\widetilde\Gamma_{\tau,k}$
defined as in~(1.12)),
which proves the property~{\it i)} of Claim~2. 
Concerning the
property~{\it ii)}, observe first that by property~{\it i)} a trajectory $\gamma_y(\cdot)$ of
$g_{\tau,k}$ starting at a point
$y\in Q_{\tau,k}\doteq\displaystyle{\overline{\widetilde\Gamma_{\tau,k}}\setminus\bigcup_{j\in\I_\tau} B_j}$
cannot esacape from~$Q_{\tau,k}$
through a point of $\partial^-\widetilde\Gamma_{\tau,k}$.
Thus, since~(3.45) shows that $|g_{\tau,k}|$ is bounded away from zero, 
and because by (3.34) one has $\partial\,Q_{\tau,k}\setminus
\partial^-\widetilde\Gamma_{\tau,k}\subset \displaystyle{\bigcup_{j\in\I_\tau} S_j}=\Sigma_\tau$,
it follows that
 $\gamma_y(\cdot)$
must cross the level set~$\Sigma_\tau$ in finite time~$\T_{\tau,k}(y)>0$,
and hence (3.63), (3.64) are verified.
In fact, with the same arguments above one can show that
every trajectory $\gamma_y(\cdot)$ starting at a point of
%
$$
Q^h_{\tau,k}\doteq\displaystyle{\overline{\widetilde\Gamma_{\tau,k}}
\setminus\bigcup_{j\in\I^h_\tau} B_j}\,,
\qquad\quad 1\leq h\leq \eta_\tau\,,
\tag 3.68
$$
%
crosses the set~$\Sigma_\tau^h\subset \Sigma_\tau$ in finite time~$\T^h_{\tau,k}(y)\geq \T_{\tau,k}(y)$. 
Next, observe that setting
$$
I_{\tau,k}\doteq
\big\{i\in\{1,\dots,N\}
\,;\ 
\P_i\cap\widetilde\Gamma_{\tau,k}\neq\emptyset\,\big\}\,,
\tag 3.69
$$
by definition (3.58) for every $i\in I_{\tau,k}\setminus(\J_{\tau,k}\cup\{k\})$
there will be some point $\overline x_i\in 
\overline{\Gamma_{\tau,k}}$
such that
$$
\big\langle\nabla V_k(\overline x_i)-\nabla V_i(\overline x_i),\, 
f(\overline x_i, v^k)\big\rangle\geq 0\,.
\tag 3.70
$$
Thus, relying on (3.43), (3.46), (3.70), we derive
$$
\aligned
\big\langle\nabla V_i(\overline x_i),\, 
f(\overline x_i, v^k)\big\rangle
&=
\big\langle\nabla V_k(\overline x_i),\, 
f(\overline x_i, v^k)\big\rangle -
\big\langle\nabla V_k(\overline x_i)-\nabla V_i(\overline x_i),\, 
f(\overline x_i, v^k)\big\rangle
\\
\noalign{\smallskip}
&< -1+ 2\ve_1\,.
\endaligned
\tag 3.71
$$
Then, since (3.27), (3.46) imply
$\overline{\widetilde\Gamma_{\tau,k}}\subset \overline{\Gamma_{\tau,k}}
\subset B(y_k, 2\ve_2)\cap B_L$,
using (3.9), (3.11),  (3.71), we find
$$
\aligned
\big\langle\nabla V_i(x),\, 
f(x, v^k)\big\rangle
&\leq
\big\langle\nabla V_i(\overline x_i),\, 
f(\overline x_i, v^k)\big\rangle +
\big|\big\langle\nabla V_i(x),\, 
f(x, v^k)\big\rangle-
\big\langle\nabla V_i(\overline x_i),\, 
f(\overline x_i, v^k)\big\rangle
\big|
\\
\noalign{\smallskip}
&< -1+ 2\ve_1+c_5 4\ve_2
\\
&< -1+ 3\ve_1
\qquad\quad\forall~x\in\overline{\widetilde\Gamma_{\tau,k}}\,,
\quad i\in I_{\tau,k}\,.
\endaligned
\tag 3.72
$$
%
%
Hence, setting  
$$
x(t)\doteq x(t;\,y, v^k)\,, 
\qquad\qquad
y\in Q^h_{\tau,k}\,,\qquad 1\leq h\leq \eta_\tau\,,
\tag 3.73
$$
and observing that,
for every fixed $0\leq s<t\leq \T^h_{\tau,k}(y)$, 
by (3.16), (3.69) there will be some index~$i(s)\in I_{\tau,k}$ such that $\tV(x(s))=V_{i(s)}(x(s))$,
relying on (3.43), (3.46), (3.64), (3.72), we derive
$$
\aligned
\tV\big(x(t)\big)&\leq V_{i(s)}\big(x(t)\big)
\\
\noalign{\smallskip}
&=V_{i(s)}\big(x(s)\big)+\int_s^t
\big\langle\nabla V_{i(s)}(x(\sigma)),\, f(x(\sigma), v^k)\big\rangle~d\sigma
\\
\noalign{\medskip}
&\leq V_{i(s)}\big(x(s)\big)+(-1+3\ve_1)(t-s)
\\
\noalign{\medskip}
&=
\tV\big(x(s)\big)+(-1+3\ve_1)(t-s)\,,
\endaligned
\tag 3.74
$$
which yields (3.65) since $\T^h_{\tau,k}(y)\geq \T_{\tau,k}(y)$.
Observe now that, by the observations at
point~{\bf 7}, we can apply Lemma~4 
for every collection of sets
$\big\{\Sigma_{\tau,k}\,;~~k\in\I^h_\tau\big\}$,
and hyperplanes $\big\{\pi_{k,i}\,;~~k,i\in\I^h_\tau\big\}$,
$h=1,\dots,\eta_\tau$,
$\tau\in\I_1^\ast$. Thus we deduce that there exists
some constant $c_{10}>0$ such that
$$
d_{\Sigma^h_\tau} \big(\partial^-\G^h_\tau\big)
\geq c_{10}\qquad\quad\forall~h=1,\dots,\eta_\tau, \quad\ \tau\in\I_1^\ast\,.
\tag 3.75
$$
Since $x\big(\T^h_{\tau,k}(y)\big)\in\Sigma^h_\tau$,
relying on (3.64), (3.75) we find that, for every fixed $\tau\in\I_1^\ast$,
  $1\leq h\leq \eta_\tau$, $k\in\I^h_\tau$,
 using the same notation in~(3.73) one has
$$
\big|y-x\big(\T^h_{\tau,k}(y)\big)\big|\geq
d_{\Sigma^h_\tau} \big(\partial^-\G^h_\tau\big)
\geq c_{10}
\qquad\quad\forall~y\in\partial^-\G^h_\tau\cap Q^h_{\tau,k}\,.
\tag 3.76
$$
On the other hand, by (3.27), (3.48), (3.64), we derive
$$
\big|y-x\big(\T^h_{\tau,k}(y)\big)\big|\leq c_9\cdot \T^h_{\tau,k}(y)
\qquad\quad\forall~y\in\partial^-\G^h_\tau\cap Q^h_{\tau,k}\,,
\tag 3.77
$$
which, together with (3.76), yields
$$
\T^h_{\tau,k}(y)\geq\frac{c_{10}}{c_9}
\qquad\quad\forall~y\in\partial^-\G^h_\tau\cap Q^h_{\tau,k}\,.
\tag 3.78
$$
Therefore, observing that by  (3.23) one has
$$
\tV\big(x(\T^h_{\tau,k}(y))\big)=\tau
\qquad\quad\forall~y\in\partial^-\G^h_\tau\cap Q^h_{\tau,k}\,,
$$
thanks to (3.78), and relying on (3.3), (3.74), we deduce
that, for every fixed $\tau\in\I_1^\ast$,
  $1\leq h\leq \eta_\tau$, $k\in\I^h_\tau$, there holds
$$
\aligned
\tV(y)&\geq \tV\big(x(\T^h_{\tau,k}(y))\big)+(1-3 \ve_1)\,\T^h_{\tau,k}(y)
\\
\noalign{\medskip}
&\geq \tau +\frac{\T^h_{\tau,k}(y)}{4}
\\
\noalign{\medskip}
&\geq \tau +\frac{c_{10}}{4c_9}
\qquad\quad\forall~y\in\partial^-\G^h_\tau\cap Q^h_{\tau,k}\,.
\endaligned
\tag 3.79
$$
Hence, since by definitions (2.40), (3.60), (3.68), one has
$$
\partial^-\G^h_\tau= \bigcup_{k\in\I^h_\tau} \partial^-\G^h_\tau\cap Q^h_{\tau,k}\,,
$$
it follows from (3.79) that
$$
\tV(y)>\tau +\ve_4
\qquad\quad\forall~y\in\partial^-\G^h_\tau\,,
\quad 1 \leq h \leq \eta_\tau\,,
\quad \tau\in\I_1^\ast\,,
\tag 3.80
$$
where $\ve_4\doteq \frac{c_{10}}{8c_9}$.
Moreover, with the same computations in (3.79) we derive also the estimates
$$
\tV(y)
>\tau
\qquad\quad\forall~y\in \G^h_\tau\setminus\bigcup_{j\in\I^h_\tau} B_j\,, 
\quad 1 \leq h \leq \eta_\tau\,,
\ \ \tau\in[\tau_0, T''[\,,
\tag 3.81
$$
$$
\tV(y)>\tau +\frac{\displaystyle{\min_h\, d_{\Sigma^h_\tau} \big(\partial^-\G^h_\tau\big)}}{8c_9}
\qquad\quad\forall~y\in\partial^-\G^h_\tau\,,
\quad 1 \leq h \leq \eta_\tau\,,
\quad \tau\in[\tau_0, T''[\,.
\tag 3.82
$$
Notice that (3.81), in particular, implies 
\, $\partial^{\strut -}\G^h_\tau\cap\Sigma^h_\tau=\emptyset$\,,\,
for all $1\leq h \leq \eta_\tau$, and hence one has
$$
\chi_\tau\doteq \tau+\frac{\displaystyle{\min_h\, d_{\Sigma^h_\tau} \big(\partial^-\G^h_\tau\big)}}{8c_9}
>\tau\,.
\tag 3.83
$$
%
To conclude, observe that by construction, for every given $\tau\in[\tau_0, T''[$,
$1 \leq h \leq \eta_\tau$,
the set
$$
\big\{x\in\R^n~;~~\tV(x)\geq\tau\big\}\setminus\partial^-\G^h_\tau
\tag 3.84
$$
consists of two connected components, one of which, say $\O^h$, contains $\Sigma_\tau^{\strut h}$.
Thus, since (3.61) implies~$\overset\circ\to\G_\tau=\displaystyle{\bigcup_{h=1}^{\eta_\tau} \overset\circ\to\G^h_\tau}$,
and because of (3.81), there holds
$$
\bigcup_{h=1}^{\eta_\tau} \O_h=\Sigma_\tau\cup
\Big(\displaystyle{\overset\circ\to\G_\tau\setminus \bigcup_{j\in\I_\tau}B_j}
\Big)\,.
\tag 3.85
$$
On the other hand, (3.80), (3.82),   imply 
$$
\gathered
\big\{x\in\R^n~;~~\tau\leq \tV(x)\leq\chi_\tau\big\}
\subset \bigcup_{h=1}^{\eta_\tau} \O_h
\qquad\quad\forall~\tau\in[\tau_0, T''[\,,
\\
\noalign{\medskip}
\big\{x\in\R^n~;~~\tau\leq \tV(x)\leq\tau+\ve_4\big\}
\subset \bigcup_{h=1}^{\eta_\tau} \O_h
\qquad\quad\forall~\tau\in\I_1^\ast\,,
\endgathered
\tag 3.86
$$
and hence (3.85), (3.86) together yield
$$
\gather
\big\{x\in\R^n~;~~\tau\leq \tV(x)\leq\chi_\tau\big\}
\subset \Sigma_\tau\cup
\Big(\displaystyle{\overset\circ\to\G_\tau\setminus \bigcup_{j\in\I_\tau}B_j}\Big)
\qquad\quad\forall~\tau\in[\tau_0, T''[\,,
\tag 3.87
\\
\noalign{\medskip}
\big\{x\in\R^n~;~~\tau\leq \tV(x)\leq\tau+\ve_4\big\}\subset
\Sigma_\tau\cup
\Big(\displaystyle{\overset\circ\to\G_\tau\setminus \bigcup_{j\in\I_\tau}B_j}\Big)
\qquad\quad\forall~\tau\in\I_1^\ast\,.
\tag 3.88
\endgather
$$
Recalling the definition (3.23) of $\Sigma_\tau$, we recover from (3.87), (3.88)
the inclusions
$$
\gather
\overset\circ\to\G_\tau\setminus\bigcup_{j\in\I_\tau} B_j\,\supset
\big\{x\in\R^n~;~\tau<\tV(x)\leq\chi_\tau\big\}
\qquad\quad\forall~\tau\in[\tau_0, T''[\,,
\tag 3.89
\\
\noalign{\medskip}
\overset\circ\to\G_\tau\setminus\bigcup_{j\in\I_\tau} 
B_j\,\supset
\big\{x\in\R^n~;~\tau<\tV(x)\leq\tau+\ve_4\big\}
\qquad\quad\forall~\tau\in\I_1^\ast\,.
\tag 3.90
\endgather
$$
which, in turn, together with (3.83), yield (3.66), (3.67), and thus we complete the proof 
of the claim.

\v

Notice that, by definitions (3.38), (3.58), (3.59), and from the above proof of Claim~2
it follows that
the inclusion in (3.90) is verified
also for all time $\tau$ in the set
$$
\gathered
\I_2^\ast\doteq
\Big\{\tau\in[\tau_0, T''[~~\setminus~ \I_1^\ast~;\ \ \ \eta_t=\eta_\tau\,, \ \ \I^h_t\subset
\I^h_\tau
\quad\
\forall~t\in [\tau,\, \tau^\sharp[\,,\ \ h=1,\dots, \eta_\tau\,
 \Big\}\,.
\endgathered
\tag 3.91
$$
Hence, we derive
$$
\tau^\sharp>\tau+\ve_4
\qquad\quad\forall~\tau\in\I_2^\ast\,.
\tag 3.92
$$
\v

\n{\bf 9.}
Relying on the properties {\it i), ii)} stated in Claim~2, for every fixed
$\tau\in[\tau_0, T''[$ we shall construct now
a patchy feedback on the open region
$$
\Omega_\tau\doteq
\displaystyle{\overset\circ\to \G_\tau\setminus \bigcup_{j\in\I_\tau} 
B_j}\,.
\tag 3.93
$$
To this end we first need to slightly enlarge some of the domains defined in (3.59).
Namely, for every~$k\in\I_\tau$, consider the set
$$
\widehat\J_{\tau,k}
\doteq
\big\{
i\in\J_{\tau,k}\cap\I_\tau\,;\ i>k\,,\ \ k\in\J_{\tau,i}
\,\big\}\,,
\tag 3.94
$$
fix some positive constant $\rho\ll\ve_3$,
denote by $\pi_{k,i}^\rho$ the hyperplane parallel
to $\pi_{k,i}$ that lies in the half space $\pi_{k,i}^+=\{ x\in\R^n\,;~~V_k(x)>V_i(x)\}$ at a distance $\rho$
from~$\pi_{k,i}$, and call $\pi_{k,i}^{\rho,-}$ 
the half space determined by $\pi_{k,i}^\rho$ that contains $\pi_{k,i}$.
Then, set
$$
\gather
\widehat \Gamma_{\tau,k}\doteq
\cases
\Gamma_{\tau,k}\cap
\displaystyle{\bigcap_{i\in\J_{\tau,k}\setminus\widehat\J_{\tau,k}} \pi_{k,i}^-}
\cap
\displaystyle{\bigcap_{i\in\widehat\J_{\tau,k}} \big(\pi_{k,i}^{\rho,-}\cap\Gamma_{\tau,i}\big)
}\quad &\text{if}\qquad \J_{\tau,k}\neq\widehat\J_{\tau,k},\ \widehat\J_{\tau,k}\neq\emptyset\,,
\\
\noalign{\smallskip}
\Gamma_{\tau,k}\cap
\displaystyle{\bigcap_{i\in\widehat\J_{\tau,k}} \big(\pi_{k,i}^{\rho,-}\cap\Gamma_{\tau,i}\big)
}\quad &\text{if}\qquad \J_{\tau,k}=\widehat\J_{\tau,k}\neq\emptyset\,,
\\
\noalign{\smallskip}\widetilde\Gamma_{\tau,k}
\quad &\text{if}\qquad \widehat\J_{\tau,k}=\emptyset\,,
\endcases
\tag 3.95
\\
\noalign{\bigskip}
\Omega_{\tau,k}\doteq 
\widehat \Gamma_{\tau,k}\setminus \bigcup_{j\in\I_\tau} B_j\,,
\tag 3.96
\endgather
$$
and observe that, by definitions (3.59), (3.62), (3.94), (3.95), (3.96), one has
$$
\gathered
\partial\,\widehat\Gamma_{\tau,k}
\setminus \bigcup
\Sb
h\in\I_{\tau}
\\
h>k \ 
\endSb\widehat\Gamma_{\tau,h}\subset
\partial\,\widetilde\Gamma_{\tau,k}\,,
\\
\noalign{\medskip}
\partial\,\Omega_{\tau,k}\setminus \Big(\Sigma_{\tau,k}\cup\bigcup
\Sb
h\in\I_{\tau}
\\
h>k \ 
\endSb\Omega_{\tau,h}
\Big)\subset\partial^-\,\widetilde\Gamma_{\tau,k}\,.
\endgathered
$$
Thus, by property {\it i)} of Claim~2 it follows that
the vector field $g_{\tau,k}(x)=f(x, v^k)$ satisfies
the inward-pointing condition~(1.5) at every point 
$x\in\partial\,\Omega_{\tau,k}\setminus \displaystyle{\Big(\Sigma_{\tau}\cup\bigcup
\Sb
h\in\I_{\tau}
\\
h>k \ 
\endSb\Omega_{\tau,h}
\Big)}$. Then, letting $g_\tau$ denote the vector field on $\Omega_\tau$
defined by
$$
g_\tau(x)\doteq g_{\tau,k}(x)
\qquad\text{if}\qquad x\in \Delta_{\tau,k}\doteq\Omega_{\tau,k}\setminus
\bigcup
\Sb
h\in\I_{\tau}
\\
h>k \ 
\endSb
\Omega_{\tau,h}\,,
\tag 3.97
$$
and considering the map $U_\tau : \Omega_\tau \to \bU$ defined by
$$
U_\tau(x)\doteq v^k\qquad\text{if}\qquad x\in \Delta_{\tau,k}\,,
\tag 3.98
$$
in view of Remark~1.3 we deduce that the triple 
$\big(\Omega_\tau,\, g_\tau,\, (\Omega_{\tau,k},\, g_{\tau,k})_{k\in\I_\tau}\big)$
is a patchy vector field on~$\Omega_\tau$ associated to the patchy feedback
$\big(\Omega_\tau,\, U_\tau,\, (\Omega_{\tau,k},\, v^k)_{k\in\I_\tau}\big)$.
%
Notice that, by definitions (3.59), (3.62), (3.94), (3.95), (3.96), (3.97), one has
$$
\Delta_{\tau,k}\subset \displaystyle{\overline{\widetilde\Gamma_{\tau,k}}\setminus\bigcup_{j\in\I_\tau} B_j}
\qquad\quad\forall~k\in\I_\tau\,,
$$
and hence we may apply the property ii) of Claim~2 to a trajectory of $g_\tau$
passing through the domain~$\Delta_{\tau,k}$.
\v

\n
{\bf Claim~3.} {\it 
The patchy vector field $g_\tau$ on the domain $\Omega_\tau$,
$\tau\in[\tau_0, T''[$\,, defined in (3.97)
enjoys the following properties.
\item{i)} 
For any $y\in \Omega_\tau$, and for every Carath\'eodory trajectory $\gamma_y(\cdot)$
of 
$$
\dot x = g_\tau (x)
\tag 3.99
$$
 starting at $y$, there exists
a time $\T_\tau(y,\,\gamma_y)>0$
so that one has
$$
\gamma_y\big(\T_\tau(y,\,\gamma_y)\big)\in\Sigma_\tau\,,
\tag 3.100
$$
and there holds
$$
t+\tV\big(\gamma_y(t)\big)\leq 
\tV(y)+3\ve_1\cdot t
\qquad\quad\forall~0\leq t\leq \T_\tau(y,\,\gamma_y)\,.
\tag 3.101
$$
\item{ii)} For any $\tau\in[\tau_0, T''[$\,, one has
$$
\tau^\sharp\doteq
\sup \Big\{t\in[\tau,\,T''[~;\quad \Sigma_s\subset \Omega_\tau\quad\forall~s\in\,]\tau,\,t]\Big\}>\tau\,.
\tag 3.102
$$
Moreover, there exists a positive constant $\ve_4$ so that
there holds
$$
\tau^\sharp>\tau+\ve_4
\qquad\quad\forall~\tau\in\I_1^\ast\cup\I_2^\ast\,.
\tag 3.103
$$
}

\v

\n{\bf 10.}\ {\smc Proof of Claim~3.} 
Given $y\in \Omega_\tau$, let $\gamma_y$ be a trajectory  of (3.99) starting at $y$,
and set 
$$
t_{\max}\big(\gamma_y\big)\doteq
\sup\big\{t>0~;~~\gamma_y \ \ \text{is \ defined \ on} \ \ [0, t]
\big\}\,.
\tag 3.104
$$
By the properties of the patchy vector fields 
recalled in Section~1 and relying on  Claim~2 one can recursively construct
two increasing sequences of times $0=t_0<t_1<\cdots<t_{\overline \nu}\leq t_{\max}$,
and of indices $i_1<i_2<\cdots<i_{\overline \nu}\in \I_\tau$
with the
following properties:
\v

\item{a)} $\gamma_y$ is a solution of $\dot x = g_{\tau, i_\nu}(x)$
taking values in $\Delta_{\tau,i_\nu}$ 
for all $t\in\,]t_{\nu-1},\,t_\nu]$, $1\leq \nu \leq {\overline \nu}$;

\item{b)} $\gamma_y(t_\nu)\in\partial\,\Omega_{\tau,i_{\nu+1}}$ for all 
$1\leq \nu < {\overline \nu}$\,, \, and $\gamma_y(t_{\overline \nu})\in\Sigma_\tau\cup \displaystyle{\bigcup
\Sb
i\in\I_{\tau}
\\
i>i_{\overline \nu}
\endSb
\partial\,\Omega_{\tau,i}}$\,;

\item{c)} $t_\nu-t_{\nu-1}<\T_{\tau,i_\nu}\big(\gamma_y(t_{\nu-1})\big)$
for all $1\leq \nu < {\overline \nu}$\,, \, and 
$t_{\overline \nu}-t_{{\overline \nu}-1}\leq\T_{\tau,i_{\overline \nu}}\big(\gamma_y(t_{{\overline \nu}-1})\big)$\,.
\v

\n
Notice that, since $\{i_\nu\}_\nu$ is strictly increasing, 
and because ${\overline \nu}\leq |\I_\tau| \leq N$
($N$ being the number of quadratic finction~$V_i$ that appear
in the definition (1.13) of the map $\tV$), we can produce
a sequence of times $t_\nu$, and of indices $i_\nu\in I_\tau$,
$1\leq \nu\leq \widehat \nu$, of such type
 so that $t_{\widehat\nu}= t_{\max}$.
Hence, since
$
\gamma_y(t_{\widehat\nu})\in\displaystyle{\bigcup
\Sb
i\in\I_{\tau}
\\
i>i_{\widehat\nu}
\endSb
\partial\,\Omega_{\tau,i}}
$
would imply that
the trajectory $\gamma_y$ could be prolonged after time
$t_{\widehat\nu}$, which is in 
 contrast with 
the maximality  of $t_{\widehat\nu}$, 
by property b) it follows that $\gamma_y(t_{\widehat\nu})\in\Sigma_\tau$,
proving (3.100). Next, applying repeatedly the estimate (3.65) of Claim~2, 
and recalling that $\gamma_y(0)=y$,
we
derive
$$
\aligned
\tV\big(\gamma_y(t)\big)&\leq \tV\big(\gamma_y(t_{\nu})\big)
+(-1+3\ve_1)\cdot \big(t-t_\nu\big)
\\
\noalign{\medskip}
&\leq\tV(y)
+(-1+3\ve_1)\cdot t
\qquad\quad\forall~t\in\,]t_{\nu-1}, t_\nu]\,,
\ \ 0< \nu\leq\widehat\nu\,,
\endaligned
$$
which yields (3.101). To conclude the proof of the claim, we only need to  observe
that, by definition~(3.93), the estimates (3.102), (3.103) are precisely the same as 
the estimates 
(3.66), (3.67), (3.92)
established at point {\bf 8}.
\v

\n{\bf 11.} Relying on Claim~3, we shall construct now a patchy feedback on the region $\D$
defined in (3.26). To this end, proceeding by induction on $m\geq 0$,
we introduce  a sequence of times $\tau_m$ defined as follows.
Observe that, by definition (3.91), for every $\tau\in[\tau_0, T''[~~\setminus~ 
\big(\I_1^\ast\cup\I_2^\ast\big)$ one has
$$
\Theta_\tau\doteq
\big\{
t\in\,]\tau,\,\tau^\sharp[~;\quad \text{either}\quad \eta_t<\eta_\tau\,,\quad\
\text{or}\quad\ |\I_t^h|>|\I^h_{\tau}| \ \ \ \text{for some} \
\ \  1\leq h\leq \eta_\tau
\big\}\neq \emptyset\,.
$$
Then, 
letting $\tau_0$ be the constant defined in (3.5), for every $m>0$,  
set
$$
\tau_m\doteq
\cases
\tau_{m-1}+\ve_4
\qquad\text{if}\quad \tau_{m-1}\in\I_1^\ast\cup\I_2^\ast\,,
\\
\noalign{\smallskip}
\inf\, \Theta_{\tau_{m-1}}
\qquad\text{otherwise\,.}
\endcases
\tag 3.105
$$
By construction, and because of (3.102), (3.103), there holds
$$
\Omega_{\tau_m}\supset
\big\{x\in\R^n~;~\tau_m<\tV(x)\leq\tau_{m+1}\big\}
\qquad\quad\forall~m\geq 0\,.
\tag 3.106
$$
Moreover, observing that $t \mapsto \eta_t$ is a decreasing map
and that $\eta_t\leq N$, $|\I_t^h|\leq N$, for all $t$ and $h$, \linebreak
it follows that  $\{\tau_m\}_{m\geq 0}$ is a strictly increasing sequence 
enjoing the property
$$
\tau_m\notin\I_1^\ast\cup\I_2^\ast
\quad
\Longrightarrow\quad
\exists~p>m\,, \quad p< m+N^2\quad\ \text{s.t.}\quad\ \tau_p\in \I_1^\ast\cup\I_2^\ast\,.
\tag 3.107
$$
In turn, (3.105), (3.107) imply that for every $m$ there exists some $p>m$,
$p< m+N^2$, such that $\tau_p>\tau_m+\ve_4$.
Thus, we deduce that there will be some integer $\mu$ such that
$\tau_{\mu}\leq T''<\tau_{\mu+1}$, and hence,
by (3.26), (3.106)
one has
$$
\D\subset\Omega\doteq\displaystyle{\bigcup_{m=0}^{\mu} \Omega_{\tau_m}}\,.
\tag 3.108
$$
Let's introduce the total ordering
$$
(m,k)\prec (p,h)\qquad\quad \text{if \ either}\qquad m>p
\qquad\text{or \ else}\qquad m=p,\quad k <h\,,
\tag 3.109
$$
on the index set
$$
\A= \big\{(m,k)~:~ m=0, \dots , \overline m,\quad k\in\I_{\tau_m}\big\}\,.
$$
Then, if we define the vector field $g$ on $\Omega$ by setting
$$
g(x)\doteq g_{\tau_m,k}(x)
\qquad\text{if}\qquad x\in D_{m,k}\doteq\Omega_{\tau_m,k}\setminus
\bigcup
\Sb
(m,k)\prec (p,h)
\endSb
\Omega_{\tau_p,h}\,,
\tag 3.110
$$
and consider the map $U : \Omega \to \bU$ defined by
$$
U(x)\doteq v^k\qquad\text{if}\qquad x\in D_{m,k}\,,
\tag 3.111
$$
in view of the observations at point {\bf 9} we deduce that the triple 
$\big(\Omega,\, g,\, (\Omega_{\tau_m,k},\, g_{\tau_m,k})_{_{(m.k)\in\A}}\big)$
is a patchy vector field on $\Omega$ associated to the patchy feedback
$\big(\Omega,\, U,\, (\Omega_{\tau_m,k},\, v^k)_{_{(m.k)\in\A}}\big)$,
so that one has
$$
g(x)=f\big(x, U(x)\big)
\qquad\quad \forall~x\in\Omega\,.
\tag 3.112
$$
Given $y\in \Omega$, let $\gamma_y$ be a Carath\'eodory trajectory  of (1.9) starting at $y$,
and define $t_{\max}\big(\gamma_y\big)$ as in~(3.104).
By the properties of the patchy vector fields 
and relying on  Claim~3 one can recursively construct
an increasing sequences of times $0=t_0<t_1<\cdots<t_{\overline \nu}\leq t_{\max}$,
and a decreasing sequence of indices $m_1>m_2>\cdots> m_{\overline \nu}$,
so that, setting $\gamma_\nu\doteq \gamma\restriction_{]t_{\nu-1},\,t_\nu]}$,
$1\leq \nu\leq \overline\nu$,
there holds:
\v

\item{a)} $\gamma_y$ is a solution of $\dot x = g_{\tau_{m_\nu}}(x)$
taking values in $\Omega_{m_\nu}$ 
for all $t\in\,]t_{\nu-1},\,t_\nu]$, $1\leq \nu \leq {\overline \nu}$;

\item{b)} $\gamma_y(t_{\nu-1})\in
\partial\,\Omega_{\tau_{m_{\nu}}}$ for all 
$1< \nu \leq {\overline \nu}$\,, \, and $\gamma_y(t_{\overline \nu})\in\Sigma_{\tau_0}\cup \displaystyle{\bigcup
\Sb
m_{\overline \nu}<p
\endSb
\partial\,\Omega_{\tau_p}}$\,;

\item{c)} $t_\nu-t_{\nu-1}<\T_{\tau_{m_\nu}}\big(\gamma_y(t_{\nu-1}), \gamma_{\nu}\big)$
for all $1\leq \nu < \overline \nu$\,, \, and 
$t_{\overline \nu}-t_{\overline \nu-1}\leq\T_{\tau_{m_{\overline \nu}}}
\big(\gamma_y(t_{\overline \nu-1}), \gamma_{\overline\nu}\big)$\,.
\v

\n
Notice that, since $\{m_\nu\}_\nu$ is strictly decreasing, 
and because ${\overline \nu}\leq \mu$, we can produce
a sequence of times $t_\nu$, and of indices $m_\nu$,
$1\leq \nu\leq \widehat \nu$, of such type
 so that $t_{\widehat\nu}= t_{\max}$.
Thus, since
$
\gamma_y(t_{\widehat\nu})\in\displaystyle{\bigcup
\Sb
m_{\overline \nu<p}
\endSb
\partial\,\Omega_{\tau_p}}
$
would imply that
the trajectory $\gamma_y$ could be prolonged after time~$t_{\widehat\nu}$, which is in 
 contrast with 
the maximality  of $t_{\widehat\nu}$, 
by property b) it follows that $\gamma_y(t_{\widehat\nu})\in\Sigma_{\tau_0}$,
and hence, by (3.29), one has $\gamma_y(t_{\widehat\nu})\in B_\ve$.
Next, given $y\in \D$,
applying repeatedly the estimate (3.101) of Claim~3, 
we
derive
$$
\aligned
\tV\big(\gamma_y(t_{\nu})\big)
&\leq \tV\big(\gamma_y(t_{\nu-1})\big)
+(-1+3\ve_1)\cdot\big(t_{\nu}-t_{\nu-1}\big)
\\
\noalign{\medskip}
&\leq\tV(y)
+(-1+3\ve_1)\cdot t_{\nu}
\qquad\quad\forall~0< \nu\leq \widehat \nu\,.
\endaligned
\tag 3.113
$$
Relying on the estimate (3.113)
in the case $\nu=\widehat\nu$, and thanks to (3.3), (3.11), (3.14), (3.26), (3.27),  
we find
$$
\aligned
t_{\widehat \nu}&\leq\frac{\tV(y)}{1-3\ve_1}
\\
\noalign{\medskip}
&\leq \big(1+2\ve_1\big)\big(V(y)+\ve_2^2\big)
\\
\noalign{\medskip}
&\leq V(y)+2\ve_1 T'+\big(1+2\ve_1\big)\ve_2^2
\\
\noalign{\medskip}
&< V(y)+\ve\,,
\endaligned
\tag 3.114
$$
which establish the conclusion of the theorem
observing that $\gamma_y$ reaches the ball $B_\ve$
within a time~$\leq t_{\widehat\nu}$
since $\gamma_y(t_{\widehat\nu})\in B_\ve$.
\fine

\vsk

\n{\medbf 4 - Appendix}
\v 
\v

We provide here a proof of the two technical lemmas stated in Section~2,
concerning the properties of lens-shaped domains of the form~(2.32)
constructed
around a collection of spheres with uniformly bounded (from above
and from below) radii.

\v
\n{\bf Proof of Lemma~3.} \ Fix $r_0>0$, and observe that the unit 
normal to  a sphere $S$ with radius $r\geq r_0$
is Lipschitz continuous with Lipschitz constant $1/r_0$:
$$
\big|\bfn_{S}(y_1)-\bfn_{S}(y_2)\big|=\frac{|y_1-y_2|}{r}\leq \frac{|y_1-y_2|}{r_0}
\qquad\forall~y_1,\,y_2\in S\,.
$$
Hence, by (2.31), and thanks to
the Lipschitz continuity of the field $g$ and of the unit 
normal $\bfn_{S}$, we deduce that there exist $\overline\ve>0$
sufficiently small, and $\overline c'>0$,  depending only on $r_0, c_0$,   and on
$\text{Lip}(g)$, so that
$$
\big\langle\bfn_{S}(x),\, g(x)\big\rangle\leq - \overline c'\,,
\qquad\quad\forall~x\in\partial^+\Gamma_{\overline\ve}\,,
\tag 4.1
$$
proving
the transversality property of  the vector field $g$ 
to the boundary $\partial^+\Gamma_{\overline\ve}$.
Next, observe that the set $\partial^-\Gamma_{\overline\ve}$ in (2.33)
is a piecewise smooth hypersurface
parametrized by 
$$
y\mapsto \Phi(y)\doteq x\big({d_{\!_{C}}}^{\!\!2}(y)-\overline\ve^2,\, y \big)\,,
\qquad\quad y\in B(C,\,\overline \ve)\cap S\,.
$$
Hence, the tangent space to $\partial^-\Gamma_{\overline\ve}$ at every
regular point $x=\Phi(y)$ of $\partial^-\Gamma_{\overline\ve}$
is the image of
the tangent space to $S$ at $y$
under 
the differential of $\Phi$,
i.e. there holds
$$
T_{\partial^-\Gamma_{\overline\ve}}(\Phi(y))=d\Phi(y)\cdot T_S(y)\,.
\tag 4.2
$$
By standard differentiability properties of the trajectories
of $\dot x = g(x)$, one finds that at the points in which
${d_{\!_{C}}}^{\!\!2}(y)$ is differentiable there holds
$$
d\Phi(y)=
\big\langle\nabla {d_{\!_{C}}}^{\!\!2}(y)
,\,\cdot\big\rangle\,g\big(\Phi(y)\big)+
X\big(
(d_{\!_{C}}(y))^2-\overline\ve^2
\big)
$$
where $X(t)$ denotes the fundamental matrix 
solution of the linear
problem $\dot v =Dg\big(x(t,y)\big)\cdot v$, 
that coincides with the identity  matrix $Id$ at
time $t=0$.
Thus, observing that at the points where
${d_{\!_{C}}}^{\!\!2}(y)$ is differentiable
one has $|\nabla {d_{\!_{C}}}^{\!\!2}(y)|\leq 2{d_{\!_{C}}}(y)$,
we obtain
$$
\aligned
\big|
d\Phi(y)-Id
\big|
&\leq 2{d_{\!_{C}}}(y)\cdot\|g\|_{\L^\infty}+
((\overline\ve^2-{d_{\!_{C}}}^{\!\!2}(y))\cdot\text{Lip}(g))\,
e^{((\overline\ve^2-{d_{\!_{C}}}^{\!\!2}(y))\cdot\text{Lip}(g))}
\\
&\leq 2\overline\ve\cdot\|g\|_{\L^\infty}+(\overline\ve^2\cdot\text{Lip}(g))\,
e^{(\overline\ve^2\cdot\text{Lip}(g))}\,.
\endaligned
\tag 4.3
$$
In turn, 
(4.3) together with (4.2) implies 
$$
\big|\bfn_{\partial^-\Gamma_{\overline\ve}}(x)-\bfn_{S}(\Phi^{-1}(x))\big|
\leq c_{11}\,\overline\ve
\tag 4.4
$$
($\bfn_{\partial^-\Gamma_{\overline\ve}}(x)$ denoting the unit normal to $\partial^-\Gamma_{\overline\ve}$),
for some constant $c_{11}>0$
depending only on 
$\|g\|_{\L^\infty},\, \text{Lip}(g)$.
Then, by the Lipschitz continuity of~$g$
and of the unit normal $\bfn_{S}$, we deduce from 
(2.31), (4.4) that, choosing
$\overline\ve>0$
sufficiently small, 
there exist some constant $\overline c''>0$, depending only on $r_0, c_0$,   and on
$\|g\|_{\L^\infty}$, $\text{Lip}(g)$,
so that at every regular point $x=\Phi(y)$ of $\partial^-\Gamma_{\overline\ve}$
there holds
$$
\big\langle\bfn_{\partial^-\Gamma_{\overline\ve}}(x),\,g(x)\big\rangle\leq - \overline c''\,.
$$
Clearly, by continuity this implies that $g(x)\in\overset\circ\to T_{\Gamma_{\overline\ve}}(x)$
at every irregular point of $\partial^-\Gamma_{\overline\ve}$ ( $\overset\circ\to T_{\Gamma_{\overline\ve}}$
denoting the interior of the tangent cone to $\Gamma_{\overline\ve}$
defined as in (1.12)),
thus showing that 
the vector field $g$ 
is inward-pointing on
the boundary $\partial^-\Gamma_{\overline\ve}$, which
completes the proof of the lemma.
\fine

\vs
\n{\bf Remark 4.1.}
Relying on the proof of Lemma~3 one can show that 
there exists some constant $c_{12}>0$
$($depending only on $r_0, \overline c$,
$\|g\|_{\L^\infty}$, 
and on $\text{Lip}(g))$, so that there holds
$$
d_S\big(x\big({d_{\!_{C}}}^{\!\!2}(y)-{\overline\ve}^2,\,y\big)\big)
> c_{12}\, \overline\ve^2
\qquad\ \ \forall~y\in B(C,\overline\ve/2)\cap S\,.
\tag 4.5
$$
Indeed, 
notice that
thanks to the Lipschitz continuity of the field $g$
we may choose the constants $\overline c',\, \overline\ve$
so that 
the estimate in (4.1) holds for all points
$x\in \Gamma_{\overline\ve}\,,$ i.e. such that
$$
\big|\langle\bfn_{S}(y),\,g(x)\rangle\big|
\geq \overline c'
\qquad\quad\forall~x\in\Gamma_{\overline\ve}\,,
\quad\forall~y\in\partial^+\Gamma_{\overline\ve}\,.
\tag 4.6
$$
Relying on (4.6) we then deduce that
$$
\aligned
d_S\big(\Phi(y)\big)&\geq
\big|\langle\Phi(y)-y,\,\bfn_{S}(y)\rangle\big|
\\
\noalign{\medskip}
&\geq \overline c '
\big(\overline\ve^2-{d_{\!_{C}}}^{\!\!2}(y)\big)
\\
\noalign{\medskip}
&>
\frac{\overline c' \ \overline\ve^2}{2}
\qquad\qquad \forall~y\in 
B(C,\overline\ve/2)\cap S\,,
\endaligned
\tag 4.7
$$
which proves (4.5), with $c_{12}\doteq\overline c'/2$,
\vs

\n{\bf Proof of Lemma~4.} \ 
\v

\n{\bf 1.} \ 
We will provide a proof of a more general result than the one stated in the lemma.
Namely, we will show that there exist constants $\overline\ve',\, c_4>0$,
 so that, 
for every given set of indices
$\I\subset\{1,\dots,\nu\}$, if we consider
the sets
$$
C^\I\doteq\bigcup_{k\in\I}
C_k\,,
\qquad\quad
\G^\I \doteq\bigcup_{k\in\I}
\overline{\Gamma_{_{\overline\ve'}}^{^{\J_k}}}\,,
\tag 4.8
$$
$$
\gathered
\partial^- \G^\I
\doteq
\partial \,\G^\I
\setminus\bigcup_{k=1}^\nu B_k\,,
\endgathered
\tag 4.9
$$
one has
$$
d_{C^\I} \big(\partial^-\G^\I\big)
\geq c_4\,.
\tag 4.10
$$
Clearly, in the particular case where $\I=\{1,\dots,\nu\}$,
we have
$$
C^\I=C\,,\qquad\quad 
\G^\I=\G\,,\qquad\quad
\partial^- \G^\I=\partial^- \G\,,
$$
and hence we recover the estimate (2.41) from (4.10).
The proof of (4.10), for an arbitrary \linebreak set~$\I\subset\{1,\dots,\nu\}$, will be obtained
proceeding by induction
on the number~$|\Pi|$ of hyperplanes contained in the set $\Pi$
considered in (2.37).
Notice that, setting
$$
\partial^{-\I}\, \Gamma_{_{\overline\ve'}}^{^{\J_k}}
\doteq
\partial\, \Gamma_{_{\overline\ve'}}^{^{\J_k}}
\setminus \bigg(
\bigcup_{j=1}^\nu B_j\cup
\bigcup
\Sb
j\in\I
\\
j\neq k\ 
\endSb
\overline{\Gamma_{_{\overline\ve'}}^{^{\J_j}}}
\bigg)
\qquad\quad \forall~k=1,\dots,\nu\,,
\tag 4.11
$$
by definitions (4.8), (4.9), one has
$$
\partial^-\G^\I =\bigcup_{k\in\I}
\overline{\partial^{-\I}\, \Gamma_{_{\overline\ve'}}^{^{\J_k}}}\,,
$$
and hence there holds
$$
d_{C^\I} \big(\partial^-\G^\I\big)\geq
\min_{k\in\I}\,
d_{C^\I}\big(\,\overline{\partial^{-\I}\, \Gamma_{_{\overline\ve'}}^{^{\J_k}}}\,\big)\,.
\tag 4.12
$$
Thus, in order to establish~(4.10), it will be sufficient to
prove by induction
on~$|\Pi|$ that there exist some constants $\overline\ve', c_4>0$,
so that
there holds
$$
d_{C^\I}\big(\partial^{-\I}\,\Gamma_{_{\overline\ve'}}^{^{\J_k}}\,\big)
> c_4
\qquad\quad\forall~k\in\I\,.
\tag 4.13
$$
\v

\n{\bf 2.} \ Consider first the case where  $\Pi=\emptyset$, i.e. assume that $\Pi_k=\emptyset$
for all $k$, and fix some set of indices~$\I\subset\{1,\dots,\nu\}$.
Then,  recalling the definitions (2.32), (2.33),
and observing that 
$$
\partial^- \Gamma_{\overline\ve}(C_k)=
\partial\, \Gamma_{\overline\ve}(C_k)\setminus {B_k}
\qquad\quad\forall~k=1,\dots,\nu\,,
\tag 4.14
$$
by (2.38), (4.11),  we have
$$
\gathered
\Gamma_{_{\overline\ve}}^{^{\J_k}}=
\Gamma_{\overline\ve}(C_k)\,,
\\
\noalign{\medskip}
\partial^{-\I}\, \Gamma_{_{\overline\ve'}}^{^{\J_k}}=
\partial^- \Gamma_{\overline\ve}(C_k)\setminus
\bigg(
\bigcup_{j\neq k} B_j\cup\bigcup
\Sb
j\in\I
\\
j\neq k\ 
\endSb
\overline{\Gamma_{\overline\ve}(C_j)}
\bigg)\,,
\endgathered
\qquad\quad\forall~k=1,\dots,\nu\,.
\tag 4.15
$$
Let $\overline\ve, c_{12}$ be the constants (depending only on $r_0$ and $g_1,\dots,g_\nu$)
provided by Lemma~3 and Remark~4.1 for all sets~$C_1,\dots,C_\nu$,
and observe that, choosing~$\overline\ve$ sufficiently small
so that
$$
\|g_k\|_{\strut\L^\infty}\, \overline\ve^2<\overline\ve/4
\qquad\quad\forall~k=1,\dots,\nu\,,
\tag 4.16
$$
and setting
$$
R_{\overline\ve}\doteq
\big\{y\in {S_k}\,;~~~~\overline\ve/2\leq d_{C_k}(y)\leq \overline \ve
\big\}\,,
$$
by the definition (2.36) of $C_k$
there holds
$$
B(R_{\overline\ve},\,  \|g_k\|_{\strut\L^\infty}\, \overline\ve^2\,)
\subset \bigcup_{j=1}^\nu\overset\circ\to B_j
\qquad\quad\forall~k=1,\dots,\nu\,.
\tag 4.17
$$
Moreover, since the solution $\tau\mapsto x(\tau, y)$ of the Cauchy problem 
$\dot x = g_k(x),$ $x(0)=y$, satisfies 
$$
|x(\tau, y)-y|\leq  \|g_k\|_{\strut\L^\infty}\,\tau
\qquad\quad\forall~\tau\,,
$$
%
%
%
we deduce from (4.17) that
$$
\gathered
\overline{\Gamma_{\overline\ve}({C_k})}\setminus
\bigcup_{j=1}^\nu\overset\circ\to B_j
\subset
\big\{\, x(\tau,y)\,; ~~y\in B({C_k},\,\overline \ve/2)\cap {S_k}\,,~~
{d_{\!_{{C_k}}}}^{\!\!\!\!\!\!2} \,(y)-\overline\ve^2
\leq \tau\leq 0\,\big\}
\cap B\big(C_k,\,2\,\overline\ve\,\big)
\\
\forall~k=1,\dots,\nu\,.
\endgathered
\tag 4.18
$$
On the other hand, by (4.14), (4.15) one has 
$$
\partial^{-\I}\,\Gamma_{_{\overline\ve}}^{^{\J_k}}\subset
\overline{\Gamma_{\overline\ve}({C_k})}\setminus
\bigcup_{j=1}^\nu\overset\circ\to B_j
\qquad\quad\forall~k=1,\dots,\nu\,.
\tag 4.19
$$
Thus, relying on (4.18), (4.19) we find
$$
\partial^{-\I}\,\Gamma_{_{\overline\ve}}^{^{\J_k}}\subset
\big\{
x\big({d_{\!_{{C_k}}}}^{\!\!\!\!\!\!2} \,(y)-{\overline\ve}^2,\,y\big)\,,~~
y\in B({C_k},\,\overline \ve/2)\cap {S_k}
\big\}
\cap\bigg(B\big(C_k,\,2\,\overline\ve\,\big)\setminus
\bigcup_{j=1}^\nu\overset\circ\to B_j\bigg)\,,
\tag 4.20
$$
which, in turn,  applying (4.5),
yields 
$$
d_{S_k}\big(\partial^{-\I}\, \Gamma_{_{\overline\ve}}^{^{\J_k}}\,\big)
>c_{12}\, {\overline\ve}^{\,2}
\qquad\quad\forall~k=1,\dots,\nu\,.
\tag 4.21
$$
Observe now that, since the radii of $S_i$, $i=1,\dots,\nu,$ are uniformly bounded by~$r_0'$,
and because the definitions (2.36),~(2.39) imply
$$
C=\partial\bigg(
\bigcup_{j=1}^\nu B_j
\bigg)\,,
$$
it follows that there will be some constants $c_{13}, c_{14}>0$, depending only on $r_0'$, such that
there holds
$$
d_{C}(y)\geq c_{13}\,
{d_{\!_{{S_k}}}}^{\!\!\!\!\!2} \,(y)
\qquad\quad \forall~y\in B\big(C_k,\,c_{14})
\setminus\Big(\bigcup_{j=1}^\nu \overset\circ\to B_j\Big)\,,
\quad \forall~k=1,\dots,\nu\,.
\tag 4.22
$$
Hence, thanks to (4.20), (4.22), choosing $\overline\ve$ sufficiently small so that 
$$
\overline\ve<\frac{c_{14}}{2}\,,
\tag 4.23
$$
and observing that
$$
d_{C^\I}(y)\geq 
d_{C}(y)\qquad\quad \forall~y\in \R^n\,,
\tag 4.24
$$
we recover from (4.21)
the estimates  (4.13),
with $c_4=c_{12}^2\cdot c_{13}\,\overline\ve^2,$ and
$\overline\ve'=\overline\ve$ satisfying~(4.16), (4.23).
\v

\n{\bf 3.} \  Given $p\geq 1$, 
suppose now that there exists some constants $\overline c_p>0$
so that, letting~$\overline\ve$ be the constant provided by Lemma~3
and satisfying (4.16), (4.23),
when $|\Pi|<p$ for every
set of indices~$\I\subset\{1,\dots,\nu\}$ there holds
$$
\align
d_{C^\I} \big(\partial^-\G^\I\big)
&\geq \overline c_p\,,
\tag 4.25
\\
\noalign{\medskip}
d_{C^\I}\big(\partial^{-\I}\,\Gamma_{_{\overline\ve}}^{^{\J_k}}\,\big)
&> \overline c_p
\qquad\quad\forall~k\in\I\,.
\tag 4.26
\endalign
$$
%
Then, consider the case where~$|\Pi|=p$. 
Fix $\I\subset\{1,\dots,\nu\}$, $k\in\I$.
Our goal is to show that there exists some constant
$\overline c_{p+1}>0$ so that
the estimate in (4.26) is verified 
with $\overline c_{p+1}$
in place of $\overline c_p$.
Clearly, if $\Pi_k=\emptyset$ we recover the estimates in  (4.26)
from the proof derived at point~{\bf 2}. Hence, we need to consider
only the case where $|\Pi_k|=|\J_k|>0$.
Then, recalling the definitions (2.32), (2.33),  by (2.38), (4.11), (4.14),
one has
$$
\Gamma_{_{\overline\ve}}^{^{\J_k}}=
\Gamma_{\overline\ve}({C_k})\cap
\bigcap_{i\in\J_k}\pi_{k,i}^-\,,
\tag 4.27
$$
%
%
$$
\gathered
\\
\noalign{\medskip}
\partial^{-\I}\,\Gamma_{_{\overline\ve}}^{^{\J_k}}
= E_1^\I \cup E_2^\I\,,
\\
\noalign{\medskip}
E_1^\I\doteq
\bigg(
\partial^- \Gamma_{\overline\ve}({C_k})\setminus
\bigg(
\bigcup_{j\neq k} B_j\cup\bigcup
\Sb
j\in\I
\\
j\neq k\ 
\endSb
\overline{\Gamma_{_{\overline\ve}}^{^{\J_j}}}\ \bigg)
\bigg)\cap \bigcap_{i\in\J_k}\pi_{k,i}^-\,,
%
\\
\noalign{\medskip}
E_2^\I\doteq
\bigcup_{i\in\J_k}E_{2,i}^\I
\qquad\quad
E_{2,i}^\I\doteq
\bigg(
\bigg(\Gamma_{\overline\ve}({C_k})
\cap
\bigcap_{j\in\J_k}
(\pi_{k,j}^-\cup\pi_{k,j})
\bigg)\setminus
\bigg(
\bigcup_{j=1}^\nu B_j\cup
\bigcup
\Sb
j\in\I
\\
j\neq k\ 
\endSb
\overline{\Gamma_{_{\overline\ve}}^{^{\J_j}}}
\bigg)
\bigg)\cap\pi_{k,i}\,.
\endgathered
\tag 4.28
$$
%
 Observe first that,
letting~$\overline\ve$ be the constant provided by Lemma~3
and satisfying (4.16), (4.23), by the proof estabilished at point {\bf 2}
one immediately deduces the inequality
$$
d_{S_k}\big(E_1^\I\big)\geq
d_{S_k}\Big(
\partial^- \Gamma_{\overline\ve}({C_k})\setminus
\bigcup_{j\neq k} B_j
\Big)>
c_{12}\, {\overline\ve}^2\,,
\tag 4.29
$$
which, together with (4.22), (4.24), yields 
$$
d_{C^\I}\big(E_1^\I\big)>
c_{12}^2\cdot c_{13}\, {\overline\ve}^4\,.
\tag 4.30
$$
%
Hence, if $E_2^\I=\emptyset$ we recover from (4.30)
 the estimates in 
(4.26) with $\overline c_{p+1}\doteq c_{12}^2\cdot c_{13}\, {\overline\ve}^4$ in place of~$\overline c_p$.
On the other hand, observe that if we let
$S_j^p$, $j\in\I$, denote the surfaces of the balls
$$
B_j^p\doteq B\big(B_j,\,{\overline c_p}/{2}\big)
\qquad j\in\I\,,
\tag 4.31
$$
%
and we consider the set
$$
{C}^{p,\strut\I}\doteq\bigcup_{k\in\I}C_k^p\,,
\qquad\qquad
C_k^p\doteq S_k^p\setminus\bigg(
\bigcup_{j\notin\I} \overset\circ\to B_j
\cup
\bigcup_{j\in\I} \overset\circ\to {B_j^p}\bigg)\,,
$$
by construction one has
$$
\partial\bigg(
\bigcup_{j\notin\I} B_j
\cup
\bigcup_{j\in\I} {B_j^p}
\bigg)=
\bigg(C\setminus\bigcup_{j\in\I} B_j^p\bigg)\cup
{C}^{p,\strut\I}\,,
\tag 4.32
$$
\v
$$
\aligned
d_{C^\I}
\Big(C\setminus\bigcup_{j\in\I} B_j^p\Big)&\geq
d_{C^\I}
\big(C\cap{C}^{p,\strut\I}\,\big)
\\
\noalign{\vskip -7pt}
&\geq d_{C^\I}
\big(\,{C}^{p,\strut\I}\,\big)
\\
&\geq\frac {\overline c_p}{2}\,.
\endaligned
\tag 4.33
$$
Thus, for every $i\in\J_k$ for which there holds
$$
\pi_{k,i}\cap
\bigg(
\bigcup_{j\in\I} B_j^p\setminus
\bigcup_{j=1}^\nu \overset\circ\to B_j
\bigg)=\emptyset\,,
\tag 4.34
$$
since the definition (4.28) implies
$$
E_{2,i}^\I\subset\bigg(
\R^n\setminus\bigcup_{j=1}^\nu B_j
\bigg)\cap \,\pi_{k,i}\,,
$$
it follows  that
$$
E_{2,i}^\I\subset\R^n\setminus\bigg(
\bigcup_{j\notin\I} B_j
\cup
\bigcup_{j\in\I} {B_j^p}
\bigg)\,.
\tag 4.35
$$
Then, relying on (4.32)-(4.33), (4.35), we deduce that
$$
\aligned
d_{C^\I}
\big(E_{2,i}^\I\big)&\geq
\min\bigg\{
d_{C^\I}\bigg(C\setminus\bigcup_{j\in\I} B_j^p\bigg),\
d_{C^\I}
\big(\,{C}^{p,\strut\I}\,\big)
\bigg\}
\\
&\geq \frac {\overline c_p}{2}\,,
\endaligned
\tag 4.36
$$
for all $i\in\J_k$ that satisfy (4.34).
Hence, in the case where (4.34) holds for all $i\in\J_k$, 
we obtain from (4.30), (4.36) the estimate in (4.26)
with 
$\overline c_{p+1}\doteq \min\{c_{12}^2\cdot c_{13}\,{\overline \ve}^4,\,c_p/2\}$ in place of $\overline c_p$.
Therefore, to complete the proof of the lemma it remains to derive
an  estimate of $d_{C^\I}\big(E_2^\I\big)$
when $E_2^\I\neq\emptyset$ and~(4.34) does not hold
for some $i\in\J_k$.
\v

\n{\bf 4.} \  With the same definitions and notations introduced at point {\bf 3},
consider a set of indices $\I\subset\{1,\dots,\nu\}$ for which
$E_2^\I\neq\emptyset$, 
and such that (4.34) is not satisfied for some $i\in\J_k$.
Set
$$
\widetilde\J_k\doteq
\bigg\{
i\in\J_k~;~~
\pi_{k,i}\cap
\bigg(
\bigcup_{j\in\I} B_j^p\setminus
\bigcup_{j=1}^\nu \overset\circ\to B_j
\bigg)\neq\emptyset
\bigg\}\,,
\tag 4.37
$$
and observe that, by the proof derived at point~{\bf 3},
there holds
$$
d_{C^\I}
\big(E_{2,i}^\I\big)\geq
\frac {\overline c_p}{2}
\qquad\quad\forall~i\in\J_k\setminus\widetilde\J_k\,.
\tag 4.38
$$
On the other hand, for every fixed
$i\in\widetilde\J_k$, by definition (4.31)
there will be
some constant \linebreak $\rho\in [0,\, (c_p/2)]$
such that, letting
$\widetilde S_j$, $j\in\I$, denote the surfaces of the balls
$$
\widetilde B_j\doteq B\big(B_j,\,\rho \big)
\qquad j\in\I\,,
\tag 4.39
$$
and considering the set
$$
{\widetilde C}^{\strut\I}\doteq\bigcup_{h\in\I}{\widetilde C}_h\,,
\qquad\qquad
{\widetilde C}_h\doteq \widetilde S_h\setminus\bigg(
\bigcup_{j\notin\I} \overset\circ\to B_j
\cup
\bigcup_{j\in\I} \overset\circ\to {\widetilde B_j}\bigg)\,,
\tag 4.40
$$
there holds
$$
\pi_{k,i}\cap {\widetilde C}^{\strut\I}\neq\emptyset\,.
$$
Then, as a first step towards
an estimate of $d_{C^\I}\big(E_{2,i}^\I\big)$
we will 
show that, setting 
$$
\vU_i\doteq {\widetilde C}^{\strut\I}\cap \pi_{k,i}\,,
\tag 4.41
$$
there holds
$$
d_{\vU_i}
\big(E_{2,i}^\I\big)\geq
\frac {\overline c_p}{2}
\qquad\quad\forall~i\in\widetilde\J_k\,.
\tag 4.42
$$
Recalling that 
$\pi_{k,i}=\pi_{i,k}$, set
$$
\gathered
\Pi^\ast\doteq
\Pi\setminus\{\pi_{k,i}\}\,,
\qquad\quad
\Pi_k^\ast\doteq
\Pi_k\setminus\{\pi_{k,i}\}\,,
\qquad\quad
\Pi_i^\ast\doteq
\Pi_i\setminus\{\pi_{k,i}\}\,,
\endgathered
\tag 4.43
$$
and 
observe that, by the properties of~$\pi_{k,i}$, one has~
$\pi_{k,i}\notin\Pi_j$ for all $j\neq k,i$,
and hence there holds $$\Pi^\ast=\bigcup_{j\neq k,i} \Pi_j\cup
\Pi_k^\ast\cup\Pi_i^\ast\,.$$ Moreover,
because of (4.43),
one has $|\Pi^\ast|<|\Pi|=p$.
Therefore, 
setting
$$
\J_i^\ast\doteq\J_i\setminus\{k\}\,,
\qquad\quad
\Gamma_{_{\overline\ve}}^{^{\J_i,\ast}}\doteq
\Gamma_{\overline\ve}({C_i})\cap
\bigcap_{j\in\J_i^\ast}\pi_{i,j}^-\,,
\qquad\quad
\I^\ast\doteq\I\setminus\{k\}\,,
\tag 4.44
$$
and defining
$$
\gathered
\G^{\ast,\I^\ast}
\doteq 
\cases
\overline{\Gamma_{_{\overline\ve}}^{^{\J_i,\ast}}}\cup
\displaystyle{\bigcup
\Sb
j\in\I^\ast\setminus\{i\}
\endSb
\overline{\Gamma_{_{\overline\ve}}^{^{\J_j}}}}
\quad &\text{if}\qquad i\in\I\,,
\\
\noalign{\smallskip}
\ \  \G^{\I^\ast}
\quad &\text{if}\qquad i\notin\I\,,
\endcases
\\
\noalign{\medskip}
\partial^-\G^{\ast,\I^\ast}
\doteq\partial\,\G^{\ast,\I^\ast}\setminus\bigcup_{j=1}^\nu B_j\,,
\endgathered
\tag 4.45
$$
by the inductive hypothesis we can apply the inequality (4.25)
in connection
with the set of hyperplanes $\Pi^\ast$ and hence, in particular, for the
set of indices $\I^\ast$ there holds
$$
d_{C^{\I^\ast}} \big(\partial^-\G^{\ast,\I^\ast}\big)
\geq \overline c_p\,.
\tag 4.46
$$
Relying on (4.46), and observing that by (4.44), one has
$$
\aligned
\overline{\Gamma_{_{\overline\ve'}}^{^{\J_i}}}\cap\pi_{i,k}&=
\overline{\big(\Gamma_{_{\overline\ve'}}^{^{\J_i,\ast}}\cap\pi_{i,k}^-
\big)}\cap\pi_{i,k}=
\overline {\Gamma_{_{\overline\ve'}}^{^{\J_i,\ast}}}\cap\pi_{i,k}\,,
\endaligned
$$
since $\pi_{i,k}=\pi_{k,i}$ we find that
$$
\aligned
B\big(C^{\I^\ast}\!,\, \overline c_p
\big)\cap\pi_{k,i}
&\subset 
\bigg(
\bigcup_{j=1}^\nu B_j\cup\overline{\Gamma_{_{\overline\ve}}^{^{\J_i,\ast}}}\cup
\displaystyle{\bigcup
\Sb
j\in\I^\ast
\endSb
\overline{\Gamma_{_{\overline\ve}}^{^{\J_j}}}}
\bigg)\cap\pi_{k,i}
\\
\noalign{\smallskip}
&\subset 
\bigg(
\bigcup_{j=1}^\nu B_j
\cup\bigcup_{j\in\I^\ast}\overline{\Gamma_{_{\overline\ve}}^{^{\J_j}}}
\bigg)\cap\pi_{k,i}\,.
\endaligned
\tag 4.47
$$
Moreover observe that,
letting
$$
{\widetilde C}^{\strut\I^\ast}\doteq\bigcup_{h\in\I^\ast}{\widetilde C}_h\,,
$$
by definitions (4.39)-(4.40) one has
$$
{\widetilde C}^{\strut\I^\ast}\subset
B\big(C^{\I^\ast}\!,\,\rho
\big)\,,\qquad\quad\rho<c_p/2\,,
$$
and hence there holds
$$
B\big({\widetilde C}^{\strut\I^\ast}\!,\,c_p/2
\big)\subset
B\big(C^{\I^\ast}\!,\,c_p\
\big)\,.
\tag 4.48
$$
Thus, (4.47), (4.48) together, yield
$$
B\big({\widetilde C}^{\strut\I^\ast}\!,\,c_p/2
\big)\cap\pi_{k,i}
\subset\bigg(
\bigcup_{j=1}^\nu B_j
\cup\bigcup_{j\in\I^\ast}\overline{\Gamma_{_{\overline\ve}}^{^{\J_j}}}
\bigg)\cap\pi_{k,i}\,.
\tag 4.49
$$
On the other hand, observe that,
by the properties of $\pi_{k,i}$, we have
$$
\widetilde {S_k}\cap \widetilde S_i=\widetilde {S_k}\cap \pi_{k,i}=\widetilde S_i\cap\pi_{k,i}\,,
$$
which, in turn, by definition (4.40) implies
$$
\widetilde {C_k}\cap \pi_{k,i}=\widetilde C_i\cap\pi_{k,i}\,.
\tag 4.50
$$
Thanks to (4.50), it follows from (4.49) that
$$
B\big(\vU_i,\,c_p/2
\big)\cap\pi_{k,i}
\subset\bigg(
\bigcup_{j=1}^\nu B_j
\cup\bigcup_{j\in\I^\ast}\overline{\Gamma_{_{\overline\ve}}^{^{\J_j}}}
\bigg)\cap\pi_{k,i}\,.
\tag 4.51
$$
Then, recalling the definition (4.28) of $E_{2,i}^\I$,
we deduce from (4.51) that
$$
E_{2,i}^\I\subset \R^n\setminus B\big(\vU_i,\,c_p/2
\big)\,,
$$
which clearly implies (4.42).
Finally, in order to obtain an estimate of $d_{C}(E_{2,i}^\I)$,
notice that, since the radii of $S_i$, $i=1,\dots,\nu,$ are uniformly bounded by~$r_0'$,
and thanks to the properties of the hyperplanes in $\Pi$,
there will be some constants $c_{15}, c_{16}>0$, depending only on $r_0'$, such that
there holds
$$
d_{{{S_k}}}(y)> c_{15}\,
d_{\!_{\vU_i}}^{\, 2}(y)
\qquad\quad \forall~y\in 
\Big(B\big(S_k,\,c_{16})
\cap
\pi_{k,i}\Big)\setminus \Big(\bigcup_{j=1}^\nu B_j\Big)\,,
\quad\forall~i\in\J_k\,, \ \ 
k=1,\dots,\nu.
\tag 4.52
$$
Therefore, choosing $\overline\ve$ so that 
$$
\overline\ve<\frac{c_{16}}{2}\,,
\tag 4.53
$$
and observing that by the same computations at point~{\bf 2}
one has 
$$
E_{2,i}^\I\subset \partial^{-\I}\,\Gamma_{_{\overline\ve}}^{^{\J_k}}
\subset
B\big(S_k, 2\overline\ve\big)\setminus \Big(\bigcup_{j=1}^\nu B_j\Big)\,,
$$
we deduce from (4.42), (4.52) that
$$
d_{{{S_k}}}\big(E_{2,i}^\I\big)\geq c_{15}\cdot
\frac {\overline c_p^2}{4}
\qquad\quad\forall~i\in\widetilde\J_k\,.
\tag 4.54
$$
Hence, letting~$\overline\ve$ be the constant provided by Lemma~3
and satisfying (4.16), (4.23), (4.53),
relying on (4.22), (4.24) we recover from (4.54) the estimates
$$
d_{C^\I}\big(E_{2,i}^\I\big)\geq c_{13}\cdot c_{15}^2\cdot
\frac {\overline c_p^4}{16}
\qquad\quad\forall~i\in\widetilde\J_k\,,
\tag 4.55
$$
which, together with (4.30), (4.38), yield
$$
\aligned
d_{C^\I} \big(\partial^{-\I}\,\Gamma_{_{\overline\ve}}^{^{\J_k}}\big) 
&\geq \min\big\{d_{C^\I} \big(E_1\big),\,\min_{i\in\J_k}d_{C^\I} \big(E_{2,i}^\I\big)\big\}
\\
\noalign{\smallskip}
&\geq \overline c_{p+1}\,,
\endaligned
\tag 4.56
$$
where $\overline c_{p+1}\doteq\min\{c_{12}^2\,c_{13}\,{\overline \ve}^4,\ \overline c_p/2, \ c_{13}\cdot c_{15}^2\cdot
\overline c_p^4/16\}$.
This establishes the estimate in (4.26) in the case where $E_2^\I\neq\emptyset$,
and (4.34) does not hold for some $i\in\J_k$,
and the proof of the lemma is completed.
\fine

\parindent 60pt

\vsk

\end